\documentclass[11pt]{article}
\usepackage[english]{babel}
\usepackage[cp850]{inputenc}
\usepackage[dvips]{graphicx}
\usepackage{amsmath,amsfonts,amsthm,amssymb}
\usepackage[usenames, dvipsnames]{color}
\usepackage{fancyhdr}
\usepackage{stmaryrd}
\usepackage[colorlinks=true,citecolor=red,linkcolor=blue,urlcolor=blue,pdfstartview=FitH]{hyperref}
\usepackage{dsfont}
\usepackage{xcolor}
\usepackage{epsfig}

\font\tenmath=msbm10 scaled 1200

\font\sevenmath=msbm7 scaled 1200
\font\fivemath=msbm5 scaled 1200

\newfam\mathfam \textfont\mathfam=\tenmath
\scriptfont\mathfam=\sevenmath \scriptscriptfont\mathfam=\fivemath

\def\R{{\mathbb R}}
\def\N{{\mathbb N}}
\def\E{{\mathbb E}}
\def\L{{\cal L}}
\def\V{{\cal V}}
\def\P{{\mathbb P}}
\def\Z{{\mathbb Z}}

\def\F{{\cal F}}
\def\G{{\cal G}}

\newtheorem{theo}{Theorem}[section]
\newtheorem{lem}{Lemma}[section]
\newtheorem{prop}{Proposition}[section]
\newtheorem{cor}{Corollary}[section]
\newtheorem{defi}{Definition}[section]
\newtheorem{pro}{Proof}[section]

\newfam\mathfam \textfont\mathfam=\tenmath
\scriptfont\mathfam=\sevenmath \scriptscriptfont\mathfam=\fivemath

\def \^#1{\if#1i{\accent"5E\i}\else{\accent"5E#1}\fi}

\def \a{\alpha}
\def \b{\beta}

\def \e{\varepsilon}
\def \g{\gamma}

\def \th{\theta}

\def \cqfd{\quad\Box}

\def \bs{\bigskip}
\def \ni{\noindent}
\def \Tr{\mbox{Tr}}

\title{\textbf{Stochastic Approximation with Averaging Innovation Applied to Finance}}

\author{\textsc{Sophie Laruelle} \thanks{Laboratoire de Probabilit\'es et Mod\`eles al\'eatoires, UMR~7599, UPMC, case 188, 4, pl. Jussieu, F-75252 Paris Cedex 5, France. E-mail: \texttt{sophie.laruelle@upmc.fr}}  
\and \textsc{Gilles Pag\`es} \thanks{Laboratoire de Probabilit\'es et Mod\`eles al\'eatoires, UMR~7599, UPMC, case 188, 4, pl. Jussieu, F-75252 Paris Cedex 5, France. E-mail: \texttt{gilles.pages@upmc.fr}}}

\date{}

\begin{document}

\maketitle

\vspace{-0.6cm}

\begin{center}
First Version: December, 10 2009 \\
This Version (accepted in Monte Carlo Methods and Applications): December, 4 2011 
\end{center}

\begin{abstract}
The aim of the paper is to establish a convergence theorem for multi-dimensional stochastic approximation when the ``innovations'' satisfy some ``light'' averaging properties in the presence of a pathwise Lyapunov function. These averaging assumptions allow us to unify apparently remote frameworks where the innovations are simulated (possibly deterministic like in Quasi-Monte Carlo simulation) or exogenous (like market data) with ergodic properties. We propose several fields of applications and illustrate our results on five examples mainly motivated by Finance.
\end{abstract}

\paragraph{Keywords} \textit{Stochastic Approximation, sequence with low discrepancy, quasi-Monte Carlo, $\alpha$-mixing process, G\`al-Koksma theorem, stationary process, ergodic control, two-armed bandit algorithm, calibration, optimal asset allocation, Value-at-Risk, Conditional Value-at-Risk}.

\bs \ni {\em 2010 AMS Classification:} 62L20,
secondary: 
37A25, 
60G10, 
93C40, 
62P05, 
91B32.

\section{Introduction}
\label{uneASM}
The aim of this paper is to establish a convergence theorem for multi-dimensional recursive stochastic approximation in a non-standard framework (compared to the huge literature on this field, see \cite{BMP}, \cite{Duf}, \cite{KusYin}, \cite{AndMouPri}, etc): we will significantly relax our assumption on the innovation process by only asking for some natural ``light'' ergodic or simply averaging assumptions, compensated by a reinforcement of the mean reversion assumption since we will require the existence of a {\em pathwise} Lyapunov function. We will show that this approach unifies seemingly remote settings: those where the innovations are simulated or even deterministic (quasi-Monte Carlo simulation) and those where the innovations are exogenous data (like market data). Especially in the latest case it may be not realistic to make a priori too stringent assumptions on the dynamics of such data process, like mixing or Markov. On the other hand, the pathwise Lyapunov assumption is definitely an intrinsic limitation to the kind of problem we can deal with, compared to the procedures extensively investigated in \cite{BMP} or more recently in \cite{DFMP} where innovations are Markovian and share mixing properties.

However, we provide several examples, mainly inspired by Finance, to illustrate the fact that the field of application of our framework is rather wide and can solve efficiently various kinds of problems, some of them having already been considered in the literature.

Let us be more specific: this paper presents convergence results for $\R^d$-valued stochastic approximation procedures of Robbins-Monro type (see \cite{RobMon} for the original paper), namely 
\begin{equation}\label{ASclassicASM}
	\th_{n+1}=\th_n-\g_{n+1}H(\th_n,Y_n), \quad n\geq0, \quad \th_0 \in\R^d,
\end{equation}	
where $H:\R^d\times\R^q\to\R^d$ is a Borel function, $(\g_n)_{n\geq1}$ a sequence of positive steps and  the ``innovation'' sequence $(Y_n)_{n\geq0}$ satisfies some ``elementary" averaging assumptions ($\th_0$ is assumed to be deterministic in this introduction for convenience simplicity). In fact, we will consider a slightly more general setting which includes an extra noisy term
\begin{equation}\label{ASpaperASM}
	\th_{n+1}=\th_n-\g_{n+1}\left(H(\th_n,Y_n)+\Delta M_{n+1}\right), \quad n\geq0,
\end{equation}	
where $(\Delta M_{n})_{n\geq1}$ is a sequence of  $\R^d$-valued $\F_n$-adapted martingale increments for a filtration $\F_n$.

To establish the $a.s.$ convergence of the sequence $(\th_n)_{n\geq0}$ toward its ``target'' $\th^*$ (to be specified later on), the idea is to make the assumption that the innovation sequence $(Y_n)_{n\geq0}$ satisfies an {\em averaging} property in a ``linear'' setting: typically that, for a wide enough class $\V$ of integrable functions (with respect to a probability measure $\nu$),
\vskip-0.3cm
\begin{equation}\label{averageASM}
\forall f\in\V, \quad \frac 1n \sum_{k=0}^{n-1}f(Y_k)\underset{n\rightarrow\infty}{\longrightarrow}\int_{{\R}^q} fd\nu
\end{equation}
at a common rate of convergence to be specified further on. If $\V\supset{\cal C}_b({\R}^q,\R)$, this implies
$$\frac{1}{n}\sum_{k=0}^{n-1}\delta_{Y_k}\overset{({\R}^q)}{\underset{n\rightarrow\infty}{\Longrightarrow}}\nu \quad a.s.$$
by a separability argument ($\overset{({\R}^q)}{\Longrightarrow}$ denotes the weak convergence of probability measures). Such a sequence $(Y_n)_{n\geq0}$ is often called ``stable'' in the literature, at least when it is a Markov chain. If $\V=L^1(\nu)$, the sequence $(Y_n)_{n\geq0}$ may be called in short ``ergodic'' although no true ergodic framework comes in the game at this stage. The target of our recursive procedure (\ref{ASpaperASM}) is then, as expected, a zero, if any, of the (asymptotic) {\em mean function} of the algorithm defined as
$$h(\th):=\int_{{\R}^q}H(\th,y)\nu(dy).$$
The key assumption is the existence of {\em pathwise Lyapunov function} with respect to the innovation $i.e.$ a function $L$ satisfying
$$\left\langle \nabla L(\th)\left.\right|H(\th,y)-H(\th^*,y)\right\rangle\geq 0$$
for every $\th$ and $y$. This assumption may look very stringent but in fact, it embodies standard framework of Stochastic Approximation with Markov representation of the form (\ref{ASclassicASM}) when the $(Y_n)_{n\geq0}$ is i.i.d. since, under appropriate integrability assumptions, it can be rewritten as follows in canonical form
$$
\th_{n+1}=\th_n-\g_{n+1}\left(\widetilde{H}(\th_n,Y_n)+\Delta \widetilde{M}_{n+1}\right),\quad n\geq0,\quad\th_0\in\R^d,
$$
where $\widetilde{H}(\th,\cdot)=h(\th)$ and $\Delta \widetilde{M}_{n+1}=H(\th_n,Y_{n})-h(\th_n)$, $n\geq0$. Then $(\Delta \widetilde{M}_n)_{n\geq1}$ is a sequence of $\sigma(Y_0,\ldots,Y_{n-1})$-martingale increments  (under appropriate integrability assumptions). Finally $\widetilde{H}(\th,\cdot)=h(\th)$ does not depend on $y$ so that the above notion of pathwise Lyapunov function reduces to the standard one. The above canonical form has been extensively investigated (and extended) in many textbooks on Stochastic Approximation (see \cite{BMP}, \cite{Duf}, \cite{KusCla}, \cite{KusYin}).

Our main theorem (Theorem \ref{Thm1ASM}) let us retrieve almost entirely the classical results about $L^p$-boundedness and $a.s.$ convergence of this procedure under standard Lyapunov assumption. Many extensions have been developed when $(Y_n)_{n\geq0}$ or even $(\th_n,Y_n)_{n\geq0}$ have a Markovian dynamics (see the seminal textbook \cite{BMP} and more recent contributions like \cite{DFMP} and several references therein). The main constraint induced by such an approach is that the existence as well as assumptions on the solution of the Poisson equation related to this chain are needed.

As a first field of applications, we are interested in quasi-random numbers. The original idea of replacing by uniformly distributed sequences (with low discrepancy) i.i.d. innovations in recursive stochastic approximation procedures goes back to the early 1990's in \cite{LapPagSab}, leading to ``Quasi-Stochastic Approximation'' ($QSA$, referring to $QMC$ for Quasi-Monte Carlo). The framework in \cite{LapPagSab} was purely one-dimensional whereas many numerical tests h proved the efficiency of $QSA$ in a multi-dimensional setting. The aim is to establish a convergence theorem in this higher dimensional setting under natural regularity assumptions ($i.e.$ based on Lipschitz regularity rather than finite variation in the Hardy \& Krause or in the measure sense, often encountered in the QMC world). As concerns the low discrepancy sequences, our framework is probably close to the most general one to get pointwise $a.s.$ convergence of stochastic approximation. 

As a second setting, we consider the case when $(Y_n)_{n\geq0}$ is a functional of $\alpha$-mixing process satisfying a priori no Markov assumption. These processes are stationary and dependent, so more realistic to model inputs made of real data. To describe the class of functions $\V$ we need to prove the convergence of the series of covariance coefficients of the innovations. To this end we use some results in \cite{PelUte} and the covariance inequality for $\alpha$-mixing process (see \cite{Dou}). Next with the probabilistic version of the G\`al-Koksma theorem  (see \cite{GalKok} and \cite{Ben,Ben1,BenPag}) we prove that this class is large enough ($L^{2+\delta}(\nu)\subset\V$, $\delta>0$). Finally we examine the case of homogeneous Markov chain with (unique) invariant distribution $\nu$. Several convergence results of stochastic approximation have been proved  in this setting  in \cite{BMP}, but they all rely on the existence (and some regularity properties)  of a solution to Poisson equation. To describe $\V$ we add an  assumption on the transition of the chain which allows us to prove that this class does not depend on the initial value of the chain.

Finally we propose several examples of applications illustrated with numerical experiment. They can be parted in two classes: the first one devoted to  simulated innovations ($i.e.$ Numerical Probability methods) and the second one deals with the applications involving real data. Primarily we present a simple case of calibration: the search,  for a derivative product in a financial model, of an implicit model parameter  fitting with  its market value. We implement the algorithm with both an i.i.d. sequence and a quasi-Monte Carlo  sequence to compare their respective rates of convergence. 
The second example is devoted to the recursive computation of risk measures commonly considered in energy portfolio management: the Value-at-Risk and the Conditional-Value-at-Risk. We design a stochastic gradient and a companion procedure to compute risk measures (like in~\cite{BarFriPag, Fri}) and we show that they can be successfully implemented in a $QSA$ framework. 
In the third example, we solve numerically a ``toy''  long term investment problem  leading to a  static potential minimization derived from an ergodic control problem (see \cite{LokZer}). The potential is related to the invariant measure of a diffusion so that the innovation relies on (inhomogeneous Markov) Euler schemes with decreasing step introduced in~\cite{LamPag} (see also~\cite{Lem}). These three   stochastic approximation  procedures relie on {\em simulated} data. The fourth example is the so-called two-armed bandit introduced in learning automata and mathematical psychology in the 1950's (see \cite{NarSha}). Its $a.s.$ behaviour in the i.i.d. setting has been extensively investigated in \cite{LamPagTar} and \cite{LamPag2} and then partially extended in~\cite{TarVan} to a more general ergodic framework. We show how the starting point of this extension appears as a consequence our  multiplicative case (Theorem \ref{Thm2ASM}).
The last example describes a model of asset allocation across liquidity pools fully developed in~\cite{LarLehPag} involving exogenous real market data, {\em a priori} sharing no Markov property but on which an averaging assumption seems natural (at least within a medium laps of time). 

\medskip
The paper is organized as follows: in Section \ref{deuxASM} are stated and proved the two main results: Theorem \ref{Thm1ASM} and  its counterpart Theorem \ref{Thm2ASM}, for multiplicative noise. Section \ref{troisASM} is devoted to quasi-Stochastic Approximation, $i.e.$ the case where the innovation process is an uniformly distributed deterministic sequence over $[0,1]^q$. Section \ref{quatreASM} is devoted to applications to random innovations, namely additive noise, mixing process (functionals of $\alpha$-mixing process), ergodic homogeneous Markov chain. Section \ref{cinqASM} presents five examples of applications including numerical illustrations, mostly in connection with Finance: implicit correlation search, recursive computation of VaR and CVaR, long term investment evaluation, two-armed bandit algorithm and optimal allocation problem (more developed in \cite{LarLehPag}).

\medskip\noindent {\bf Notations} $\left\langle \cdot\left|\right.\cdot\right\rangle$ denotes the canonical Euclidean inner product and $\left|\cdot\right|$ its related norm on $\R^d$. The almost sure convergence will be denoted by $\stackrel{a.s.}{\longrightarrow}$ and $\stackrel{({\R}^q)}{\Longrightarrow}$ will denote the weak convergence of probability measures on $({\R}^q,\mathcal{B}or({\R}^q))$. $\Delta a_n=a_n-a_{n-1}$ for every sequence $(a_n)_n$.

\section{Algorithm design and main theoretical result}
\label{deuxASM}

In this paper, we consider the following general framework for recursive stochastic algorithms of the following form
\begin{equation}\label{AlgoStoASM}
\th_{n+1}=\th_n-\g_{n+1}\left(H(\th_n, Y_{n})+\Delta M_{n+1}\right), \quad n\geq 0,
\end{equation}
where $(Y_n)_{n\geq0}$ is an ${\R}^q$-valued sequence of $\F_n$-adapted random variables and $(\Delta M_n)_{n\geq1}$ is a sequence of $\F_n$-adapted martingale increment, all defined on a same filtered probability space $\left(\Omega,\F,(\F_n)_{n\geq0},\P\right)$. Moreover $\th_0\in L^1_{\R^d}(\Omega,\F_0,\P)$ and $\th_0$ is independent of $(Y_n,\Delta M_{n+1})_{n\geq0}$. The positive step sequence $(\g_n)_{n\geq1}$ is non-increasing and $H$ is a Borel function from ${\R}^d\times{\R}^q$ to ${\R}^d$.

In the following, we adopt a kind of compromise by assuming that $(Y_n)_{n\geq0}$ is a process satisfying some averaging properties and that the function $H(\th^*,\cdot)$ belongs to a class of functions (to be specified further on) for which a rate of convergence ($a.s.$ and in $L^p$) holds in (\ref{averageASM}). Moreover we need to reinforce the Lyapunov condition on the pseudo-mean function $H$ which limits, at least theoretically, the range of application of the method.

\subsection{Framework and assumptions}

Let $(Y_n)_{n\geq0}$ be an ${\R}^q$-valued random variables sequence. We will say that the sequence $(Y_n)_{n\geq0}$ satisfies a {\em $\nu$-stability assumption} or equivalently is {\em $\nu$-averaging} if
\begin{equation}\label{StableASM}
\P(d\omega)\mbox{-}a.s. \quad \frac{1}{n}\sum_{k=0}^{n-1}\delta_{Y_k(\omega)}\overset{({\R}^q)}{\underset{n\rightarrow\infty}{\Longrightarrow}}\nu.
\end{equation}

We will see that the stochastic approximation procedure defined by (\ref{AlgoStoASM}) is a recursive zero search of the (asymptotic) mean function
\vskip-0.3cm
\begin{equation}
	h(\th):=\int_{{\R}^q} H(\th,y)\nu(dy).
\label{fmoyASM}
\end{equation}
\noindent Let $p\in[1,\infty)$ and let $(\e_n)_{n\geq0}$ be a sequence of nonnegative numbers such that 
\begin{equation}\label{HepsASM}
\e_n\underset{n\rightarrow\infty}{\longrightarrow}0 \quad\mbox{and}\quad \lim_n\inf n\e_n=0.
\end{equation}
We denote by $\V_{\e_n,p}$ the class of functions which convergence rate in (\ref{averageASM}) in both $a.s.$ and in $L^p(\P)$ sense is $\e_n^{-1}$, namely 
\vskip-0.3cm
\begin{equation}
	\V_{\e_n,p}=\left\{f\in L^p(\nu) \left|\right. \frac{1}{n}\sum_{k=0}^{n-1}f(Y_k)-\int fd\nu\stackrel{\P\mbox{-}a.s. \ \& \ L^p(\P)}{=}O(\e_n)\right\}.
\label{ClasseBetaPASM}
\end{equation}

\subsection{Main result}

Now we are in a position to state an $a.s.$-convergence theorem ``\`a la" Robbins-Siegmund.
\begin{theo}\label{Thm1ASM}
$(a)$ {\rm Boundedness}. Let $h: {\R}^d \rightarrow {\R}^d$ satisfying (\ref{fmoyASM}), $H: {\R}^d\times {\R}^q \rightarrow  {\R}^d$ a Borel function and let $(Y_n)_{n\geq0}$ be a $\nu$-stable sequence ($i.e.$ satisfying (\ref{StableASM})). 
Assume there exists a continuously differentiable function $L:\R^d\rightarrow\R_+$ satisfying
\begin{equation}\label{LyapunovASM}
	\nabla L \ is \ Lipschitz \ continuous \ and \ \left|\nabla L \right|^2\leq C\left(1+L\right)
\end{equation}
and that the pseudo-mean function $H$ satisfies the \textnormal{pathwise Lyapunov} assumption 
\begin{equation}\label{lmr1ASM}
\begin{array}{c}
	\forall \th\in{\R}^d\backslash\{\th^*\}, \ \forall y\in{\R}^q, \quad	\left\langle \nabla L(\th)\left|\right. H(\th,y)-H(\th^*,y)\right\rangle\geq 0.
\end{array}
\end{equation}
Let $p\in[1,\infty)$ and let   $(\e_n)_{n\ge 1}$ be a sequence satisfying~(\ref{HepsASM}). Assume that
\begin{equation}\label{GH2ASM}
	 H(\th^*,\cdot)\in\V_{\e_n,p}.
\end{equation}
Moreover, assume that $H$ satisfies the following (quasi-)linear growth assumption
\begin{equation}
	\forall \th\in{\R}^d,\; \forall y\in{\R}^q, \quad \left|H(\th,y)\right| \leq C_H \phi(y)(1+L(\th))^{\frac{1}{2}}
\label{HHASM} 
\end{equation}
and that the martingale increments sequence $(\Delta M_{n})_{n\geq1}$  satisfies for every $n\ge0$, 
\begin{equation}
\P\mbox{-}a.s.\,\left\{\begin{array}{lll}
  \E\left(\left|\Delta M_{n+1}\right|^{2\vee\frac{p}{p-1}}\left|\right.\F_n\right) \leq C_M \phi(Y_n)^{2\vee\frac{p}{p-1}}(1+L(\th_n))^{1\vee\frac{p}{2(p-1)}} & \mbox{if} & p>1, \\
  \frac{|\Delta M_{n+1}|}{(1+L(\th_n))^{\frac 12}}\leq C_M& \mbox{if} & p=1
\end{array}\right.
\label{DeltaMASM} 
\end{equation}
where $C_M$ is a positive real constant and $\sup_{n\geq0}\left\|\phi(Y_n)\right\|_{2\vee\frac{p}{p-1}}<+\infty$.

\noindent Let $\g=(\g_n)_{n\geq1}$ be a nonnegative non-increasing sequence of ``admissible" gain parameters satisfying
\begin{equation}\label{gammaASM}
	\sum_{n\geq1}\g_n=+\infty, \quad n\e_n\g_n\underset{n\rightarrow\infty}{\longrightarrow} 0, \quad \mbox{and} \quad \sum_{n\geq1}n\e_n\max\left(\g_n^2,\left|\Delta\g_{n+1}\right|\right)<+\infty.
\end{equation}
Then, the recursive procedure defined by (\ref{AlgoStoASM}) satisfies $(L(\th_n))_{n\geq0}$ is $L^1$-bounded, $L(\th_n)\underset{n\rightarrow\infty}{\longrightarrow}L_{\infty}<+\infty$ $a.s.$, $\Delta\th_n\underset{n\rightarrow\infty}{\longrightarrow}0$ $a.s.$ and
\vskip-0.3cm
$$\sum_{n\geq1}\left\langle \nabla L(\th_n)\left|\right. H(\th_n,Y_n)-H(\th^*,Y_n)\right\rangle<+\infty.$$
$(b)$ {\rm $A.s.$ convergence toward $\th^*$}. Furthermore, if $\{\th^*\}$ is a connected component of $\{L=L(\th^*)\}$ and the pseudo-mean function $H$ satisfies the \textnormal{strict pathwise Lyapunov} assumption
\begin{equation}\label{lmrASM}
\begin{array}{c}
	\forall \delta>0, \ \forall \th\in{\R}^d\backslash\{\th^*\}, \ \forall y\in{\R}^q, \quad
	\left\langle \nabla L(\th)\left|\right. H(\th,y)-H(\th^*,y)\right\rangle\geq \chi_{_\delta}(y)\Psi_{\delta}(\th)
\end{array}
\end{equation}
where $\nu(\chi_{_\delta})>0$, $ \Psi_{\delta}$ is l.s.c. and positive on ${\R}^d\backslash\{\th^*\}$ and 
$\bigcap_{\delta>0}\{\Psi_{\delta}=0\}=\{\th^*\}$, then
$$\th_n\overset{a.s.}{\underset{n\rightarrow\infty}{\longrightarrow}}\th^*.$$
\end{theo}

\noindent {\bf Remark.} The conditions on the step sequence $\g=(\g_n)_{n\geq1}$ and $(\varepsilon_n)_{n\geq1}$ are satisfied for example by
\begin{equation}\label{GamExASM}
	\e_n=n^{-\beta}, \quad \beta\in (0,1], \quad \mbox{and} \quad \g_n=\frac{c}{n^{a}}, \quad 1-\beta<a\le 1, \quad  c>0.
\end{equation}

\paragraph{Proof.}  \textbf{{\em \underline{First ste}p}}: We introduce the function
$$\Lambda(\th):=\sqrt{1+L(\th)}$$
as a Lyapunov function instead of $L(\th)$ like in the classical case. It follows from the fundamental formula of calculus that there exists $\xi_{n+1}\in\left(\th_n, \th_{n+1}\right)$ such that
\begin{eqnarray*}
\Lambda(\th_{n+1})&=&\Lambda(\th_{n})+\left\langle \nabla\Lambda(\th_{n})\left.\right|\Delta\th_{n+1}\right\rangle+\left\langle \nabla\Lambda(\xi_{n+1})-\nabla\Lambda(\th_{n})\left.\right|\Delta\th_{n+1}\right\rangle \\
				   &\leq&\Lambda(\th_{n})+\left\langle \nabla\Lambda(\th_{n})\left.\right|\Delta\th_{n+1}\right\rangle+| \nabla\Lambda(\xi_{n+1})-\nabla\Lambda(\th_{n})|  |\Delta\th_{n+1}|. 	
\end{eqnarray*}
\begin{lem}\label{Lemme1ASM} The new Lyapunov function $\Lambda$ satisfies the two following properties
\begin{itemize}
	\item[(i)] $\nabla\Lambda$ is bounded (so that $\Lambda$ is Lipschitz).
	\item[(ii)] $\forall \th,\th'\in{\R}^d, \quad \left|\nabla\Lambda(\th')-\nabla\Lambda(\th)\right|\leq C_L\frac{\left|\th'-\th\right|}{\Lambda(\th)}$.
\end{itemize}
\end{lem}

\noindent {\bf Proof of Lemma \ref{Lemme1ASM}.} $(i)$ $\nabla\Lambda=\frac{\nabla L}{2\sqrt{1+L}}$ is bounded by (\ref{LyapunovASM}), consequently $\Lambda$ is Lipschitz.

\medskip
\noindent $(ii)$  Let $\th, \th'\in{\R}^d$, 
\begin{eqnarray*}	
\left|\nabla\Lambda(\th)-\nabla\Lambda(\th')\right|&\leq&\frac{\left|\nabla L(\th)-\nabla L(\th')\right|}{2\sqrt{1+L(\th)}}+\frac{\left|\nabla L(\th')\right|}{2}\left|\frac{\sqrt{1+L(\th')}-\sqrt{1+L(\th)}}{\sqrt{1+L(\th)}\sqrt{1+L(\th')}}\right| \\
										&\leq&\frac{[\nabla L]_{{\rm Li}p}}{2\sqrt{1+L(\th)}}\left|\th-\th'\right|+\frac{C}{2\sqrt{1+L(\th)}}[\Lambda]_{{\rm Lip}}\left|\th-\th'\right|\\
										&\leq&\frac{1}{2}\left([\nabla L]_{{\rm Lip}}+C[\Lambda]_{{\rm Lip}}\right)\frac{\left|\th-\th'\right|}{\Lambda(\th)} \\
										&=&C_L\frac{\left|\th-\th'\right|}{\Lambda(\th)}.\hspace{7cm}\cqfd
\end{eqnarray*}

\noindent Thus, applying the above lemma to $\theta =\theta_n$ and $\theta'=\xi_{n+1}$, and noting that $|\xi_{n+1}-\theta_n|\le |\Delta\th_{n+1}|$ yields 
\begin{eqnarray*}
	\Lambda(\th_{n+1})&\leq&\Lambda(\th_{n})-\g_{n+1}\left\langle \nabla\Lambda(\th_{n})\left.\right|H(\th_{n},Y_{n})\right\rangle-\g_{n+1}\left\langle \nabla\Lambda(\th_{n})\left.\right|\Delta M_{n+1}\right\rangle+C_L\frac{\left|\Delta\th_{n+1}\right|^2}{\sqrt{1+L(\th_n)}}\\
					   &=&\Lambda(\th_{n})-\g_{n+1}\left\langle \nabla\Lambda(\th_{n})\left.\right|H(\th_{n},Y_{n})-H(\th^*,Y_n)\right\rangle-\g_{n+1}\left\langle \nabla\Lambda(\th_{n})\left.\right|H(\th^*,Y_n)\right\rangle \\
					   & &-\g_{n+1}\left\langle \nabla\Lambda(\th_{n})\left.\right|\Delta M_{n+1}\right\rangle+C_L\g_{n+1}^2\frac{\left|H(\th_{n},Y_{n})+\Delta M_{n+1}\right|^2}{\sqrt{1+L(\th_n)}}.
\end{eqnarray*}
\noindent We have for every $n\geq0$, 
$$\quad \left|\g_{n+1}\left\langle \nabla\Lambda(\th_{n})\left.\right|H(\th^*,Y_n)\right\rangle\right|\leq C_{\Lambda}\g_{n+1}\phi(Y_n)\in L^1(\P)$$
since $\nabla \Lambda$ is bounded. Besides $\E\left[\left\langle \nabla\Lambda(\th_{n})\left.\right|\Delta M_{n+1}\right\rangle\left|\right.{\cal F}_n\right]=0$, $n\geq0$,
since $\Delta M_n$ is a true martingale increment and $\nabla \Lambda$ is bounded. Furthermore, owing to (\ref{HHASM}) and (\ref{DeltaMASM})
$$\E\left[\frac{\left|H(\th_{n},Y_{n})+\Delta M_{n+1}\right|^2}{\sqrt{1+L(\th_n)}}\left|\right.{\cal F}_n\right]\leq C\phi^2(Y_n)\Lambda(\th_n)$$
(where conditional expectation is extended to positive random variables). Consequently,
\begin{equation}\label{CondExpectASM}	
\begin{array}{rcl}
\E\left[\Lambda(\th_{n+1})\left.\right|\F_n\right]&\leq&\Lambda(\th_{n})\left(1+C'_L\g_{n+1}^2\phi(Y_n)^2\right)-\g_{n+1}\left\langle \nabla\Lambda(\th_{n})\left.\right|H(\th_{n},Y_{n})-H(\th^{*},Y_{n})\right\rangle \\
&&-\g_{n+1}\left\langle \nabla\Lambda(\th_{n})\left.\right|H(\th^{*},Y_{n})\right\rangle.
\end{array}
\end{equation}
\noindent We set $V_n:=\displaystyle\frac{A_n}{B_n}$, $n\ge 1$ where
$$A_n:=\Lambda(\th_{n})+\sum_{k=0}^{n-1}\g_{k+1}\left\langle\nabla\Lambda(\th_{k})\left.\right|H(\th_{k},Y_{k})-H(\th^{*},Y_{k})\right\rangle, \quad B_n:=\prod_{k=1}^n\left(1+C'_L\g_{k}^2\phi(Y_{k-1})^2\right).$$
\noindent Using the mean-reverting assumption (\ref{lmrASM}) implies that $(A_n)_{n\geq0}$ is a nonnegative process and $B_n$ is $\F_{n-1}$-adapted, $n\geq1$. Elementary computations first show that 
$$\E\left[A_{n+1}\left|\right.\F_n\right]\leq A_n\frac{B_{n+1}}{B_n}-\g_{n+1} \left\langle \nabla\Lambda(\th_n)\left.\right|H(\th^{*},Y_n)\right\rangle$$
which finally yields 
\begin{equation}\label{surmgleASM}
\forall n\geq0, \quad \E\left[V_{n+1}\left.\right|\F_n\right]\leq V_n-\Delta W_{n+1},
\end{equation}
where $W_n:=\sum_{k=0}^{n-1} \widetilde{\g}_{k+1} \left\langle \nabla\Lambda(\th_{k})\left.\right|H(\th^{*},Y_{k})\right\rangle$ with $\widetilde{\g}_n:=\frac{\g_n}{B_n}, \ n\geq0$.  

\medskip
\noindent \textbf{{\em \underline{Second ste}p}}: Now our aim is to prove that the sequence $(W_n)_{n\geq0}$ is $L^1$-bounded and $a.s.$ converges. To this end we set $S^*_n:=\sum_{k=0}^{n-1} H(\th^*,Y_k)$, then it follows
$$W_n=\sum_{k=0}^{n-1} \widetilde{\g}_{k+1} \left\langle \nabla\Lambda(\th_k)\left.\right|\Delta S^*_{k+1}\right\rangle=\widetilde{\g}_n\left\langle \nabla\Lambda(\th_{n-1})\left.\right|S^*_n\right\rangle-\sum_{k=1}^{n-1}  \left\langle S^*_k\left.\right|\widetilde{\g}_{k+1}\nabla\Lambda(\th_{k})-\widetilde{\g}_{k}\nabla\Lambda(\th_{k-1})\right\rangle.$$
\noindent First, since $\nabla\Lambda$ is bounded, note that
$$\widetilde{\g}_n\left|\nabla\Lambda(\th_{n-1})\right|\left|S^*_n\right|\leq\left\|\nabla\Lambda\right\|_{\infty}n\e_n\widetilde{\g}_n\frac{\left|S^*_n\right|}{n\e_n}\leq\left\|\nabla\Lambda\right\|_{\infty}n\e_n\g_n\frac{\left|S^*_n\right|}{n\e_n}$$
which $a.s.$ goes to 0 as $n$ goes to infinity since $n\e_n\g_n\underset{n\rightarrow\infty}{\longrightarrow}0$ by (\ref{gammaASM}) and $\left(\displaystyle\frac{S^*_n}{n\e_n}\right)_{n\geq1}$ remains $a.s.$ bounded. Moreover
\vskip-0.5cm
$$\E\left[\widetilde{\g}_n\left|\nabla\Lambda(\th_{n-1})\right|\left|S^*_n\right|\right]\leq n\e_n\g_n\left\|\nabla\Lambda\right\|_{\infty}\left\|\frac{S^*_n}{n\e_n}\right\|_1$$
which converges to $0$ in $L^1$ because $n\e_n\g_n\underset{n\rightarrow\infty}{\longrightarrow}0$ and $H(\th^*,\cdot)\in\V_{\e_n,p}$. On the other hand, 
$$\sum_{k=1}^{n-1} \left\langle S^*_k\left.\right|\widetilde{\g}_{k+1}\nabla\Lambda(\th_{k})-\widetilde{\g}_{k}\nabla\Lambda(\th_{k-1})\right\rangle=\sum_{k=1}^{n-1} \left\langle S^*_k\left.\right|\nabla\Lambda(\th_{k})\right\rangle\Delta\widetilde{\g}_{k+1}+\sum_{k=1}^{n-1} \widetilde{\g}_k\left\langle S^*_k\left.\right|\nabla\Lambda(\th_{k})-\nabla\Lambda(\th_{k-1})\right\rangle.$$
\noindent As $\nabla \Lambda=\displaystyle\frac{\nabla L}{2\sqrt{1+L}}$ is bounded by construction, we have
$$\sum_{k=1}^n\left| \Delta\widetilde{\g}_{k+1}\left\langle S^*_k\left.\right|\nabla\Lambda(\th_{k})\right\rangle\right|\leq\sum_{k=1}^n \left|\Delta\widetilde{\g}_{k+1}\right|\left|S^*_k\right|\left\|\nabla\Lambda\right\|_{\infty}\leq\left\|\nabla\Lambda\right\|_{\infty}\sum_{k=1}^n k\e_k\left|\Delta\widetilde{\g}_{k+1}\right|\left|\frac{S^*_k}{k\e_k}\right|.$$
Now, using that $\frac{a}{1+a}\leq\sqrt{a}, \ a>0$,
$$\left|\Delta\widetilde{\g}_{k+1}\right|\leq\left|\Delta\g_{k+1}\right|+\g_{k}\frac{C'_L\g_{k+1}^2\phi(Y_{k})^2}{B_{k+1}}\leq\left|\Delta\g_{k+1}\right|+\g_{k}\frac{C'_L\g_{k+1}^2\phi(Y_{k})^2}{1+C'_L\g_{k+1}^2\phi(Y_{k})^2}\leq\left|\Delta\g_{k+1}\right|+C'_L\g_{k}\g_{k+1}\phi(Y_{k}).
$$
Hence
$$\sum_{k=1}^n\left| \Delta\widetilde{\g}_{k+1}\left\langle S^*_k\left.\right|\nabla\Lambda(\th_{k})\right\rangle\right|\leq\left\|\nabla\Lambda\right\|_{\infty}\left(\sum_{k=1}^n k\e_k\left|\Delta\g_{k+1}\right|\left|\frac{S^*_k}{k\e_k}\right|+C'_L\sum_{k=1}^n k\e_k\g_k\g_{k+1}\phi(Y_k)\left|\frac{S^*_k}{k\e_k}\right|\right).$$
\noindent By H\"older's Inequality
$$\E\left(\phi(Y_{k})\left|\frac{S^*_k}{k\e_k}\right|\right)\leq\left\|\phi(Y_{k})\right\|_{\frac{p}{p-1}}\left\|\frac{S^*_k}{k\e_k}\right\|_p.$$
\noindent As $\left(\displaystyle\frac{S^*_k}{k\e_k}\right)_{n\geq0}$ is bounded, $\g$ is admissible and
$\displaystyle\sup_{k\geq0}\left\|\phi(Y_k)\right\|_{\frac{p}{p-1}}<+\infty$, then the series $\sum_{k=1}^n \Delta\widetilde{\g}_k\left\langle S^*_k\left.\right|\nabla\Lambda(\th_{k})\right\rangle$ is absolutely converging in $L^1(\P)$. 

\noindent We study now the series $\displaystyle\sum_{k=1}^n \widetilde{\g}_k\left\langle S^*_k\left.\right|\nabla\Lambda(\th_{k})-\nabla\Lambda(\th_{k-1})\right\rangle$. We have
$$\left|\nabla\Lambda(\th_k)-\nabla\Lambda(\th_{k-1})\right|\leq C'_L\frac{\left|\Delta\th_k\right|}{\sqrt{1+L(\th_{k-1})}}\leq C'_L\g_k\frac{\left|H(\th_{k-1},Y_{k-1})\right|+\left|\Delta M_{k}\right|}{\sqrt{1+L(\th_{k-1})}}.$$
\noindent We are interested in the $L^1$-convergence of the series 
$$\sum_{k=1}^n \g_k^2\left|S^*_k\right|\frac{\left|H(\th_{k-1},Y_{k-1})\right|}{\sqrt{1+L(\th_{k-1})}} \quad \mbox{and} \quad \sum_{k=1}^n \g_k^2\left|S^*_k\right|\frac{\left|\Delta M_k\right|}{\sqrt{1+L(\th_{k-1})}}.$$
\noindent For the first sum, as $\displaystyle\frac{\left|H(\th_{k-1},Y_{k-1})\right|}{\sqrt{1+L(\th_{k-1})}}\leq C_H\phi(Y_{k-1})$, we then come to
$\sum_{k=1}^n C_H\g_k^2\E\left[\left|S^*_k\right|\left|\phi(Y_{k-1})\right|\right]$ and by H\"older's inequality we obtain
$$\E\left[\left|S^*_k\right|\left|\phi(Y_{k-1})\right|\right]\leq\left\|S^*_k\right\|_p\left\|\phi(Y_{k-1})\right\|_{\frac{p}{p-1}}<+\infty$$
because $\left\|S^*_k\right\|_p=O\left(k\e_k\right)$ by (\ref{ClasseBetaPASM}) and $\displaystyle\sup_{n\geq0}\left\|\phi(Y_n)\right\|_{\frac{p}{p-1}}<+\infty$. Furthermore, as $\displaystyle\sum_{k\geq1}k\e_k\g_k^2<+\infty$ by (\ref{gammaASM}), then the series $\displaystyle\sum_{k=1}^n \g_k^2\left|S^*_k\right|\frac{\left|H(\th_{k-1},Y_{k-1})\right|}{\sqrt{1+L(\th_{k-1})}}$ converges in $L^1$.

\noindent For the second sum, we derive from H\"older's inequality
(with $p$ and $\frac{p}{p-1}$) and~(\ref{DeltaMASM})  that
$$\E\left[\left|S^*_k\right|\frac{\left|\Delta M_k\right|}{\sqrt{1+L(\th_{k-1})}}\right]\leq\left\|S^*_k\right\|_p\left\|\frac{\left|\Delta M_k\right|}{\sqrt{1+L(\th_{k-1})}}\right\|_{\frac{p}{p-1}}\leq C_M\left\|S^*_k\right\|_p\left\|\phi(Y_{k-1})\right\|_{\frac{p}{p-1}\vee 2}<+\infty.
$$
This yields that $\displaystyle\sum_{k=1}^n \g_k^2\left|S^*_k\right|\frac{\left|\Delta M_k\right|}{\sqrt{1+L(\th_{k-1})}}$ converges in $L^1$ too. Finally we then obtain that $W_n\overset{a.s.}{\underset{n\rightarrow\infty}{\longrightarrow}}W_{\infty}$ and $\sup_{n\geq0}\left\|W_n\right\|_1<+\infty$. Thus we have that 
$$(V_n+W_n)^-\leq W_n^-\leq\left|W_n\right|\in L^1(\P) \quad \mbox{since} \quad\sup_{n\geq0}\left\|W_n\right\|_1<+\infty.$$ 
As $V_0=\Lambda(\th_0)\leq C(1+\left|\th_0\right|)\in L^1$, it follows by induction from (\ref{surmgleASM}) that, for every $n\geq0$, $\E\,V_n<+\infty$. Hence $S_n:=V_n+W_n$, $n\geq0$, is a true supermartingale with an $L^1$-bounded negative part. We then deduce that
$$S_n\overset{a.s.}{\underset{n\rightarrow\infty}{\longrightarrow}}S_{\infty}\in L^1.$$ 
Now $W_n\overset{a.s.}{\underset{n\rightarrow\infty}{\longrightarrow}}W_{\infty}$ implies $V_n\overset{a.s.}{\underset{n\rightarrow\infty}{\longrightarrow}}V_{\infty}<+\infty \ a.s.$  

\medskip
\noindent \underline{{\em \textbf{Third ste}p}}: Now we show that the product $B_n$ converges $a.s.$ to derive that $A_n$ converges $a.s.$. 

$$\mbox{In fact}\hskip5cm\sum_{n\geq1}\g_n^2\phi^2(Y_{n-1})<+\infty \ a.s.,\hskip6cm$$
since $\displaystyle\sup_{n\geq0}\E\left[\phi^2(Y_{n})\right]<+\infty$ and $\sum_{n\geq1}\g_n^2<+\infty$ by combining (\ref{HepsASM}) and (\ref{gammaASM}), which in turn implies that $B_n\overset{a.s.}{\underset{n\rightarrow\infty}{\longrightarrow}}B_{\infty}<+\infty$. As a consequence $A_n\overset{a.s.}{\underset{n\rightarrow\infty}{\longrightarrow}}A_{\infty}<+\infty.$ Therefore using the mean reverting property (\ref{lmrASM}) of $H$ with respect to $\nabla\Lambda$, we classically derive that
\begin{equation}\label{GH4ASM}
\sum_{n\geq1} \g_n \left\langle\nabla\Lambda(\th_{n-1})\left.\right|H(\th_{n-1},Y_{n-1})-H(\th^*,Y_{n-1})\right\rangle<+\infty \quad a.s.
\end{equation}
$$\mbox{Consequently}\hskip4cm\Lambda(\th_n)\overset{a.s.}{\underset{n\rightarrow\infty}{\longrightarrow}}\Lambda_{\infty}<+\infty \quad a.s.\hskip6cm$$ 
\noindent As $\displaystyle\lim_{\left|\th\right|\rightarrow +\infty}L(\th)=+\infty$, $\displaystyle\lim_{\left|\th\right|\rightarrow +\infty}\Lambda(\th)=+\infty$, then the sequence $(\th_n)_{n\geq0}$ is $a.s.$-bounded and
$$
L(\th_n)\overset{a.s.}{\underset{n\rightarrow\infty}{\longrightarrow}}L_{\infty}<+\infty \quad a.s.
$$ 
\noindent Now let us show that $\Delta\th_n\underset{n\rightarrow+\infty}{\longrightarrow}0$. In fact, $\left|\Delta\th_{n+1}\right|^2\leq C\g_{n+1}^2\left(\left|H(\th_n,Y_n)\right|^2+\left|\Delta M_{n+1}\right|^2\right)$, so that
$$
\E\left[\left|\Delta\th_{n+1}\right|^2\left|\right.\F_n\right]\leq C\g_{n+1}^2\phi(Y_n)^2(1+L(\th_n))
$$
and $(L(\th_n))_{n\geq0}$ being $a.s.$ bounded,
$$\sum_{n\geq0}\E\left[\left|\Delta\th_{n+1}\right|^2\left|\right.\F_n\right]<+\infty. \ a.s.
$$
which classically implies that $\sum_{n\geq0}\left|\Delta\th_{n+1}\right|^2<+\infty \ a.s.$ 

\medskip
\noindent \textbf{{\em \underline{Fourth ste}p}}: To prove the convergence of $\th_n$ toward $\th^*$, we use Assumptions (\ref{lmrASM}) and (\ref{GH4ASM}) to deduce that
\begin{equation}\label{chiBASM}
\sum_{n\geq1} \g_n \chi_{_\delta}(Y_{n-1})\Psi_{\delta}(\th_{n-1})<+\infty \quad a.s.
\end{equation}
$$\mbox{Now,}\hskip 4cm\sum_{k=0}^n \g_{k+1}\chi_{_\delta}(Y_k)=\sum_{k=0}^n \g_{k+1}\Delta S^{\chi}_k=\g_{n+1} S^{\chi}_n-\sum_{k=1}^{n-1}\Delta\g_{k+2} S^{\chi}_{k}\hskip 4cm$$
where $S^{\chi}_n=\sum_{k=0}^n \chi_{_\delta}(Y_k)$ and we set $S^{\chi}_0=0$ and $\Delta S^{\chi}_0=0$. 

\noindent By Assumption (\ref{StableASM}), $\displaystyle\frac{S^{\chi}_n}{n}\rightarrow\nu(\chi_{_\delta})>0$
as $n\rightarrow\infty$. Let $n_0$ be the smallest integer such that
$$\forall n\geq n_0, \quad \frac{S^{\chi}_n}{n}\geq\epsilon_0=\frac{\nu(\chi_{_\delta})}{2}>0.$$
\noindent Then, a standard discrete integration by part yields
$$
\forall n\geq n_0, \quad \sum_{k=n_0}^n \g_{k+1}\chi_{_\delta}(Y_k)=n\g_{n+1}\frac{S^{\chi}_n}{n}-C_{n_0}+\sum_{k=n_0}^{n-1}k(-\Delta\g_{k+2})\frac{S^{\chi}_k}{k} \quad a.s.,
$$
where $C_{n_0}=\g_{n_0+1}S^{\chi}_{n_0-1}$. Therefore, using that the sequence $(-\Delta\g_n)_{n\geq1}$ is nonnegative,
\begin{eqnarray*}
	\sum_{k=n_0}^n \g_{k+1}\chi_{_\delta}(Y_k)&\geq&n\g_{n+1}\epsilon_0-C_{n_0}+\sum_{k=n_0}^{n-1}k(-\Delta\g_{k+2})\epsilon_0=\epsilon_0\Big(n\g_{n+1}+\sum_{k=n_0}^{n-1}k(-\Delta\g_{k+2})\Big)-C_{n_0} \\
	                                      &=&\epsilon_0\Big(\g_{n+1}+n_0\g_{n_0+1}+\sum_{k=n_0+1}^{n-1}\g_{k+1}\Big)-C_{n_0} \\
\end{eqnarray*}
by a reverse discrete integration by parts. Finally	                                      
$$\sum_{k=n_0}^n \g_{k+1}\chi_{_\delta}(Y_k)\geq\epsilon_0\Big(\g_{n+1}+\sum_{k=n_0+1}^{n-1}\g_{k+1}\Big)-C_{n_0}\rightarrow\infty \quad \mbox{as} \quad n\rightarrow\infty$$
since $\displaystyle\sum_{n\geq1}\g_n=+\infty$. We have then shown that 
$$\sum_{k\geq0} \g_{k+1}\chi_{_\delta}(Y_k)=+\infty \quad a.s.$$
\noindent Combining this fact with (\ref{chiBASM}) classically implies that 
$$\lim\inf_n\widetilde{\Psi}_{\delta}(\th_n)=0.$$
\noindent Let $\Theta_{\infty}$ be the set of limiting points of the sequence $(\th_n)_{n\geq0}$. $\Theta_{\infty}$ is a compact connected set since $(\th_n)_{n\geq0}$ is bounded and $\Delta\th_n\underset{n\rightarrow\infty}{\longrightarrow}0$. Moreover $\{\Psi_{\delta}=0\}$ is closed because $\Psi_{\delta}\leq0$ and l.s.c. and $\Theta_{\infty}$ is closed too. So $\Theta_{\infty}\cap\{\Psi_{\delta}=0\}$ is a family of nonempty compact sets which decreases as $\delta\searrow0$ since it is bounded. As a consequence, 
$$\bigcap_{\delta>0}\left(\Theta_{\infty}\cap\{\Psi_{\delta}=0\}\right)\neq\emptyset.$$
The other assumption on $\Psi_{\delta}$ implies $$\bigcap_{\delta>0}\left(\Theta_{\infty}\cap\{\Psi_{\delta}=0\}\right)\subset\bigcap_{\delta>0}\{\Psi_{\delta}=0\}=\{\th^*\},$$
so that in fact it is reduced to $\th^*$. Hence $\th^*$ is a limiting point of $(\th_n)_{n\geq0}$ which implies that $L(\th_n)$ converges towards $L(\th^*)$. By the assumption on the Lyapunov function $L$, $\{\th^*\}$ is a connected component of $\{L=L(\th^*)\}$ and as $\Theta_{\infty}$ is connected, $\Theta_{\infty}=\{\th^*\}$. Therefore
$$
\hskip 6 cm \th_n\stackrel{a.s.}{\longrightarrow}\th^* \quad \mbox{as} \quad n\rightarrow\infty. \hskip 6,25 cm\hfill\cqfd
$$

\paragraph{{\sc Back to the i.i.d. innovation setting.}} $\rhd$ Theorem \ref{Thm1ASM} contains the {\em standard martingale approach} ``\`a la" Robbins-Siegmund in the i.i.d. setting. Indeed, we consider the recursive procedure
$$\th_{n+1}=\th_n-\g_{n+1}K(\th_n,Y_{n+1}), \quad n\geq0,$$
where $(Y_n)_{n\geq1}$ is i.i.d. with distribution $\nu$ and $\th_0$ is independent of $(Y_n)_{n\geq1}$ (all defined on $(\Omega, {\cal A},\P)$). We set $\F_n=\sigma(\th_0,Y_1,\ldots,Y_n)$, $n\geq0$, $p=2$,
$$H(\th,y)=h(\th)\quad\mbox{with}\quad h(\th)=\int_{{\R}^q} K(\th,y)\nu(dy)\quad\mbox{and}\quad\Delta M_{n+1}=K(\th_n,Y_{n+1})-\E\left[K(\th_n,Y_{n+1})\left|\right.\F_n\right].$$
Assume that
$$\forall\th\in{\R}^d, \quad \left\|K(\th,Y_1)\right\|_2\leq C_K(1+L(\th))^{\frac 12}.$$
Then Assumption~(\ref{HHASM}) is satisfied by $h$ and (\ref{DeltaMASM}) holds (with $\phi\equiv1$). Furthermore, by combining (\ref{HepsASM}) and (\ref{gammaASM}), we retrieve  the step assumption in the standard Robbins-Monro Theorem, namely
\begin{equation*}
\sum_{n\geq1}\g_n=+\infty\quad\mbox{and}\quad\sum_{n\geq1}\g_n^2<+\infty.
\end{equation*}
$\rhd$ Another (naive) way to apply Theorem~\ref{Thm1ASM} in this i.i.d. setting is to focus, under the above assumption on  the averaging property so that: then $H=K$ and $\Delta M_n\equiv0$.  We still consider the above procedure but we assume furthermore the existence of a  pathwise Lyapunov function. By noticing that $(K(\th^*,Y_n))_{n\geq1}$ is i.i.d. and in $L^2$, it follows from the quadratic law of large numbers (at rate $n^{-\frac 12}$) and the Law of Iterated Logarithm at rate $O(\e_n)$ with $\e_n=\sqrt{\frac{\log\log n}{n}}$ that $K(\th^*,\cdot)\in\V_{\e_n,2}$. As a consequence the condition (\ref{gammaASM}) is clearly more restrictive than the above regular one, however any step of the form $\g_n=\frac{c}{n^{\alpha}}$, $c>0$, $\frac{3}{4}<\alpha\leq1$ satisfies~(\ref{gammaASM}).

\subsection{The case of multiplicative noise}

If we assume that the function $H$ in~(\ref{AlgoStoASM}) is of the following form
\begin{equation}\label{multiplicativeASM}
\forall\th\in{\R}^d, \, \forall y\in{\R}^q, \quad H(\th,y)=\chi(y)h(\th)+H(\th^*,y),
\end{equation}
where $\chi$ is a Borel function such that $\nu(\chi)=1$, $\chi\in\V_{\e_n,p}$ and $\sup_{n\geq0}\left\|\chi(Y_n)\right\|_{2\vee\frac{p}{p-1}}<+\infty$, $H(\th^*,\cdot)\in\V_{\e_n,p}$ and $\sup_{n\geq0}\left\|H(\th^*,Y_n)\right\|_{2\vee\frac{p}{p-1}}<+\infty$, $h$ is Lipschitz bounded with $h(\th^*)=0$, then we replace the growth assumption (\ref{HHASM}) on $H$ by one on the mean function $h$, $i.e.$
\begin{equation}\label{HyphASM}
	\forall \th\in{\R}^d, \forall y\in{\R}^q, \quad \left|h(\th)\right| \leq C_h\sqrt{1+L(\th)}
\end{equation}
and the pathwise mean-reverting assumption (\ref{lmrASM}) is the classical
\begin{equation}\label{nmrASM}
\forall\th\in{\R}^d\setminus\{\th^*\}, \quad \left\langle\nabla L\left|\right.h\right\rangle(\th)>0.
\end{equation}
\begin{theo}\label{Thm2ASM}
The recursive procedure (\ref{AlgoStoASM}) with the function $H$ defined by (\ref{multiplicativeASM}) and the previous assumptions on $\chi$ and (\ref{HyphASM})-(\ref{nmrASM}) on $h$ satisfies
$$\th_n\overset{a.s.}{\underset{n\rightarrow\infty}{\longrightarrow}}\th^*.$$
\end{theo}

\paragraph{Proof.}\textbf{{\em \underline{First ste}p}}: This setting cannot be reduced to the general setting. We use the same notations as in the proof of Theorem \ref{Thm1ASM}. With the new form of the function $H$, we obtain
\begin{eqnarray*}
\Lambda(\th_{n+1})&\leq&\Lambda(\th_{n})-\g_{n+1}\left\langle \nabla\Lambda(\th_{n})\left.\right|\chi(Y_n)h(\th_n)\right\rangle-\g_{n+1}\left\langle \nabla\Lambda(\th_{n})\left.\right|H(\th^*,Y_n)\right\rangle \\
					   & &-\g_{n+1}\left\langle \nabla\Lambda(\th_{n})\left.\right|\Delta M_{n+1}\right\rangle+C_L\g_{n+1}^2\frac{\left|H(\th_{n},Y_{n})+\Delta M_{n+1}\right|^2}{\sqrt{1+L(\th_n)}}.
\end{eqnarray*}
\noindent By the same arguments as before we get
\begin{equation*}	
\E\left[\Lambda(\th_{n+1})\left.\right|\F_n\right]\leq\Lambda(\th_{n})\left(1+C'_L\g_{n+1}^2\phi(Y_n)^2\right)-\g_{n+1}\left\langle \nabla\Lambda(\th_{n})\left.\right|\chi(Y_n)h(\th_{n})\right\rangle-\g_{n+1}\left\langle \nabla\Lambda(\th_{n})\left.\right|H(\th^{*},Y_{n})\right\rangle.
\end{equation*}
\noindent We set $V_n:=\displaystyle\frac{A_n}{B_n}$, where $A_n:=\Lambda(\th_{n})+\sum_{k=0}^{n-1}\g_{k+1}\left\langle\nabla\Lambda\left.\right|h\right\rangle(\th_{k})$ and $B_n:=\prod_{k=1}^n\left(1+C'_L\g_{k}^2\phi(Y_{k-1})^2\right)$. Using the mean-reverting assumption (\ref{nmrASM}) implies that $(A_n)_{n\geq0}$ is a nonnegative process whereas $(B_n)_{n\geq0}$ is still $\F_{n-1}$-adapted. Elementary computations show that
$$\E\left[A_{n+1}\left|\right.\F_n\right]\leq A_n\frac{B_{n+1}}{B_n}-\g_{n+1} \left\langle \nabla\Lambda(\th_n)\left.\right|H(\th^{*},Y_n)\right\rangle-\g_{n+1}\widetilde{\chi}(Y_n)\left\langle \nabla\Lambda\left.\right|h\right\rangle(\th_n)$$
where $\widetilde{\chi}(Y_n):=\chi(Y_n)-\nu(\chi)$, $n\geq0$. Finally we have
\begin{equation}\label{surmgle2ASM}
\forall n\geq0, \quad \E\left[V_{n+1}\left.\right|\F_n\right]\leq V_n-\Delta W_{n+1}-\Delta Z_{n+1},
\end{equation}
where $W_n:=\sum_{k=0}^{n-1} \widetilde{\g}_{k+1} \left\langle \nabla\Lambda(\th_{k})\left.\right|H(\th^{*},Y_{k})\right\rangle$ and $Z_n:=\sum_{k=0}^{n-1} \widetilde{\g}_{k+1} \widetilde{\chi}(Y_k)\left\langle \nabla\Lambda\left.\right|h\right\rangle(\th_{k})$ 
with $\widetilde{\g}_n:=\frac{\g_n}{B_n}, \ n\geq0$. 

\medskip
\noindent \textbf{{\em \underline{Second ste}p}}: Following the lines of the proof of Theorem \ref{Thm1ASM} we show that the sequence $(W_n)_{n\geq0}$ is $L^1$-bounded and $a.s.$ converges. Now our aim is to prove the same results for the sequence $(Z_n)_{n\geq0}$. To this end we set $S^{\widetilde{\chi}}_n:=\sum_{k=0}^{n-1} \widetilde{\chi}(Y_k)$, then it follows
$$Z_n=\sum_{k=0}^{n-1} \widetilde{\g}_{k+1}\Delta S^{\widetilde{\chi}}_{k+1}\left\langle \nabla\Lambda\left.\right|h\right\rangle(\th_k)=\widetilde{\g}_n S^{\widetilde{\chi}}_{n}\left\langle \nabla\Lambda\left.\right|h\right\rangle(\th_{n-1})-\sum_{k=1}^{n-1}S^{\widetilde{\chi}}_{k}\left(\widetilde{\g}_{k+1}\left\langle \nabla\Lambda\left.\right|h\right\rangle(\th_{k})-\widetilde{\g}_{k}\left\langle\nabla\Lambda\left.\right|h\right\rangle(\th_{k-1})\right).$$
\noindent By the same methods as for the sequence $(W_n)_{n\geq0}$ ($i.e.$ using assumptions on $H$, $\Lambda$ and $(\g_n)_{n\geq1}$), we obtain that
$$Z_n\overset{a.s.}{\underset{n\rightarrow\infty}{\longrightarrow}}Z_{\infty} \quad \mbox{and} \quad \sup_{n\geq0}\left\|Z_n\right\|_1<+\infty.$$
Thus we have that 
$$(V_n+W_n+Z_n)^-\leq (W_n+Z_n)^-\leq\left|W_n+Z_n\right|\in L^1(\P) \quad \mbox{since} \quad\sup_{n\geq0}\left\|W_n+Z_n\right\|_1<+\infty.$$ 
As $V_0=\Lambda(\th_0)\leq C(1+\left|\th_0\right|)\in L^1$, it follows by induction from (\ref{surmgleASM}) that, for every $n\geq0$, $\E\,V_n<+\infty$. Hence $S_n:=V_n+W_n+Z_n$, $n\geq0$, is a true supermartingale with a $L^1$-bounded negative part. We then deduce that
$$S_n\overset{a.s.}{\underset{n\rightarrow\infty}{\longrightarrow}}S_{\infty}\in L^1.$$ 
Now $W_n\overset{a.s.}{\underset{n\rightarrow\infty}{\longrightarrow}}W_{\infty}$ and $Z_n\stackrel{a.s.}{\longrightarrow}Z_{\infty}$ imply that $V_n\overset{a.s.}{\underset{n\rightarrow\infty}{\longrightarrow}}V_{\infty}<+\infty \ a.s.$ 

\medskip
\noindent \textbf{{\em \underline{Third ste}p}}: Like in the proof of Theorem \ref{Thm1ASM}, we have that $B_n\overset{a.s.}{\underset{n\rightarrow\infty}{\longrightarrow}}B_{\infty}<+\infty$ which implies that $A_n\overset{a.s.}{\underset{n\rightarrow\infty}{\longrightarrow}}A_{\infty}<+\infty.$ Therefore using the mean-reverting property (\ref{nmrASM}) of $h$ with respect to $\nabla\Lambda$, we classically derive that
\begin{equation}\label{GH5ASM}
\sum_{n\geq0} \g_{n+1}\nu(\chi)\left\langle\nabla\Lambda\left.\right|h\right\rangle(\th_n)<+\infty \quad a.s.
\end{equation}
The end of the proof follows the lines of the one of Theorem \ref{Thm1ASM}.\hfill$\cqfd$

\section{Application to quasi-stochastic approximation}
\label{troisASM}

This section is devoted to quasi-random innovations: the innovation sequence $(Y_n)_{n\ge0}$ becomes a deterministic uniformly distributed (u.d.) sequence $(\xi_{n+1})_{n\ge0}$ over a unit hypercube $[0,1]^q$, $i.e.$
$$\th_{n+1}=\th_n-\g_{n+1}H(\th_n,\xi_{n+1}),\quad n\geq0,\quad\th_0\in\R^d.$$
We extend  the  one-dimensional result first introduced in \cite{LapPagSab}  (see more recently \cite{MehMey,ShiMey}) to a general multi-dimensional setting with unbounded function $H$. We first recall few definitions and properties of u.d. sequences  (see~\cite{Nie} and the reference therein). We emphasize how to  apply Theorem \ref{Thm1ASM} when $H$ has ``bounded variation" on $[0,1]^q$ and when $H$ is Lipschitz continuous.

\subsection{Definition and characterisation}
\begin{defi} A $[0,1]^q$-valued sequence $(\xi_n)_{n\geq1}$ is uniformly distributed (u.d.) on $[0,1]^q$ if
$$\frac{1}{n}\sum_{k=1}^n\delta_{\xi_k}\stackrel{\left({\R}^q\right)}{\Longrightarrow}{\cal U}([0,1]^q) \quad \mbox{as} \ n\rightarrow\infty.$$
\end{defi}
The proposition below provides a characterisation of uniform distribution for a sequence $(\xi_n)_{n\geq1}$.
\begin{prop} $(a)$ Let $(\xi_n)_{n\geq1}$ be a $[0,1]^q$-valued sequence. Then $(\xi_n)_{n\geq1}$ is uniformly distributed on $[0,1]^q$ if and only if 
\vskip-0.3cm
$$D^*_n(\xi):=\sup_{x\in[0,1]^q}\Big|\frac{1}{n}\sum_{k=1}^n\mathds{1}_{\llbracket 0,x\rrbracket}(\xi_k)-\prod_{i=1}^q x^i\Big|\longrightarrow0 \quad \mbox{as} \ n\rightarrow\infty,$$
where $D^*_n(\xi)$ is called the discrepancy at the origin or star discrepancy.	

\noindent $(b)$ There exists sequences, called {\em sequences with low discrepancy} such that $D^*_n(\xi) =\Big(\frac{\log^dn}{n}\Big)$.  We refer to~\cite{Nie, BoLe} for examples of such sequences (like Halton, Kakutani, Sobol' sequences, etc).
\end{prop}

\subsection{Standard classes \texorpdfstring{$\V_{\e_n,1}$}{V(e,1)}  for quasi-stochastic approximation}

We set here $Y_n=\xi_{n+1}$, $\F_n=\{\emptyset,\Omega\}$ and $\Delta M_{n+1}\equiv0$, $n\geq0$. The strong Lyapunov condition on $H$ is crucial here. Note that the function $\phi$ becomes useless since we always consider the case $p=1$. To apply Theorem~\ref{Thm1ASM}, we mainly need to specify  the accessible classes  $\V_{\e_n,1}$ in such a framework. 

\medskip 
\noindent $\rhd$ {\bf Function with finite variation.} A function $f:[0,1]^q\to \R$ has finite variation in the measure sense if there exists a signed measure $\nu$ on $([0,1]^d,{\cal B}or([0,1]^q)$ such that $\nu(\{0\})=0$ and 
\[
\forall\, x\in [0,1]^q,\quad f(x)= f(\mbox{\bf 1})+\nu([\![0,\mbox{\bf 1}-x]\!])
\]
where $[\![x,y]\!]=\prod_{i=1}^q[x^i,y^i]$ if $x\le y$ (componentwise) and is empty otherwise and $\mbox{\bf 1}=(1,\ldots,1)$. The variation $V(f)$ of $f$is then defined as $|\nu|([0,1]^q)$ where $\nu$ denotes the total variation measure attached to $\nu$. For further details on this notion of variation, see~\cite{BoLe}. When $q=1$ this notion coincide with left continuous functions with finite variations. As concerns the slightly more general  notion of finite variation in the the Hardy and Krause sense, see~\cite{Nie} and the references therein). The role of finite variation is emphasized by the following error bound.

\begin{prop}[Koksma-Hlawka Inequality] \label{KokHlaASM}Let $\xi=(\xi_1,\ldots,\xi_n)\in ([0,1]^q)^n$ and let $f$ be a function with finite variation $V(f)$, either in the Hardy \& Krause or in the measure sense. Then
$$
\Big|\frac{1}{n}\sum_{k=1}^n f(\xi_k)-\int_{[0,1]^q}f(u)\lambda_q(du)\Big|\leq V(f)D^*_n(\xi).
$$
\end{prop}
\noindent Hence, if $(\xi_n)_{n\ge1}$ has a low discrepancy, the class $\V=\left\{f:[0,1]^q\rightarrow{\R} \mbox{ s.t. } V(f)<+\infty\right\}$ of functions with finite variations  satisfies 
$\V \subset\V_{\e_n,1}$ with $\e_n=\frac{(\log n)^q}{n}$. Consequently, if $H(\th,\cdot)\in\V$, the assumptions on admissible (non-increasing) step sequences $(\gamma_n)_{n\ge 1}$  in Theorem~\ref{Thm1ASM} reads
$$
\sum_{n\geq1}\g_n=+\infty, \quad\g_n(\log n)^q\rightarrow0 \quad\mbox{and}\quad \sum_{k\geq1}\max(\left|\Delta\g_{n+1}\right|,\g_n^2)(\log n)^q<+\infty.
$$
so that the choice of $\g_n:=\frac{c}{n^{\rho}}$, $\frac{1}{2}<\rho\leq 1$,  is admissible (like in the  i.i.d. setting).

\medskip 
\noindent $\rhd$ {\bf Lipschitz continuous functions.} If $q\geq2$ it is difficult to check whether  $f$ has finite variations in any sense: in fact these functions become ``rare'' as $q$ increases. If we look for  more natural regularity assumption to be satisfied by $H(\th^*,\cdot)$ like Lipschitz continuity, the following theorem due to Proinov (see \cite{Pro}) provides an alternative (but less ``attractive") error bound. 

\begin{theo}\label{ProinovASM} (Proinov) Assume ${\R}^q$ is equipped with the $\ell^{\infty}$-norm $\left|(x^1,\ldots,x^q)\right|_{\infty}:=\max_{1\leq i\leq q} \left|x^i\right|$. Let $\left(\xi_1,\ldots,\xi_n\right)\in\left ([0,1]^q\right)^{\otimes n}$. For every continuous function $f:[0,1]^q \rightarrow {\R}$,
$$
\Big|\frac{1}{n}\sum_{k=1}^n f\left(\xi_k\right)-\int_{[0,1]^q} f(u)\lambda_q(du)\Big| \leq C_q w_f\! \Big(D^*_n\left(\xi_1,\ldots,\xi_n\right)^{\frac{1}{q}}\Big)
$$
where $\displaystyle w_f(\delta):= \sup_{x,y\in[0,1]^q,\left|x-y\right|_{\infty}\leq\delta} \left|f(x)-f(y)\right|$, $\delta\in(0, 1)$,
is the $\ell^\infty$-uniform continuity modulus of $f$ and $C_q\in(0,\infty)$ is a universal constant only depending on $q$. If $q=1$, $C _q=1$ and if $q\geq2$, $C_q\in [1, 4]$.
\end{theo} 
\noindent As a consequence ${\rm Lip}([0,1]^q,\R)\subset \V_{\e_n,1}$ with $\e_n=\frac{\log n}{n^{\frac{1}{q}}}$ (with obvious extensions to H\" older functions). Consequently, if $H(\th^*,\cdot)\in\V$, the assumptions on admissible (non-increasing) step sequences $(\gamma_n)_{n\ge 1}$  in Theorem~\ref{Thm1ASM} reads
$$
\sum_{n\geq1}\g_n=+\infty, \quad \g_n(\log n)n^{1-\frac{1}{q}}\rightarrow0 \quad \mbox{and}  \quad \sum_{k\geq1}\max(\left|\Delta\g_{n+1}\right|,\g_n^2)(\log n )n^{1-\frac{1}{q}}<+\infty.
$$	
so that the choice of $\g_n:=\frac{c}{n}$ is always admissible (more generally $\g_n = cn^{-\rho}$, $1-\frac 1q<\rho\le1$). An  application of ``quasi-Stochastic Approximation" is proposed in Section~\ref{ICQSA} (see also~\cite{Fri}).

\section{Applications to different types of random innovations}
\label{quatreASM}

This section is devoted to some first applications of the above theorem. By applications, we mean here printing out some classes of random innovation processes $(Y_n)_{n\geq0}$ for which the averaging rate assumption (\ref{ClasseBetaPASM}) is naturally satisfied by ``large" class ${\cal V}_{\e_n,p}$.

First we present a simple framework of stochastic approximation where the noise is additive which is studied in \cite{DedDou} with some mixing properties, but here we only need (\ref{StableASM}). We showed in \cite{LarLehPag} how easily our result applies to real life stochastic optimization problem (as far as convergence is concerned).

Afterwards we focus on mixing innovations: we consider that the sequence $(Y_n)_{n\geq0}$ is a functional of a stationary $\alpha$-mixing process (satisfying condition on the summability of the mixing coefficients). 

The last application is the case of an homogeneous Markov chain which can be seen as a possible more elementary counterpart of some (convergence) result obtained $e.g.$ in \cite{BMP}.  Some (quasi-optimal) $a.s.$ rate of convergence can be obtained if $H$ is smooth enough in $\th$ (see~\cite{Lar}), but to establish a regular Central Limit Theorem it is most likely that we cannot avoid to deal with the Poisson equation.

\subsection{Recursive procedure with additive noise}

We consider here the case where the function $H$ is the sum of the mean function $h$ and a noise, namely
$$\forall \th\in{\R}^d, \quad \forall y\in{\R}^q, \quad H(\th,y)=h(\th)+y, \quad \mbox{and} \quad \Delta M_{n+1}\equiv0.$$
In this framework, the Lyapunov assumption (\ref{lmrASM}) becomes classical involving only the mean function $h$, namely 
$$\forall\th\in{\R}^d\setminus\{\th^*\}\, \quad \left\langle \nabla L\left|\right.h\right\rangle(\th)>0.$$
Likewise, the growth control assumption (\ref{HHASM}) amounts to
$$\forall\th\in{\R}^d, \quad \left|h(\th)\right|\leq C_h\sqrt{1+L(\th)},$$
provided the moment assumption $\sup_n\left\|Y_n\right\|_{\frac{p}{p-1}}<+\infty$, for some $p\in(1,\infty]$, is satisfied (take $\phi(y):=\left|y\right|\vee1$). The martingale is vanishing in this example. Finally the step assumption (\ref{gammaASM}) is ruled by the averaging rate of the sequence $(Y_n)_{n\geq0}$. 

\subsection{Functional of a stationary \texorpdfstring{$\alpha$}{alpha}-mixing process}

Here we provide a short background on $\alpha$-mixing processes and their functionals. Our motivation here is to relax as much as possible our assumption on $(Y_n)_{n\geq0}$ in order to apply stochastic approximation methods to exogenous possibly non Markovian stationary data.

We aim now at applying our convergence theorem to input sequences $(Y_n)_{n\geq0}$ which are (causal) functionals of an $\alpha$-mixing process. Consider a stationary ${\R}^q$-valued process $X=(X_k)_{k\in\mathbb{Z}}$, its natural filtration $\F_n=\F^X_n:=\sigma(X_k,\,k\leq n)$ and $\G_n=\G^X_n:=\sigma(X_k,\,k\geq n)$. The $\alpha$-mixing coefficients are defined as follows
\begin{equation}\label{MixingASM}
\alpha_n=\sup\left\{\left|\P(U\cap V)-\P(U)\P(V)\right| ,\; U\in\F_{k}, V\in\G_{k+n}, k\geq0\right\}.
\end{equation}
Let $f$ be a measurable mapping from $({\R}^q)^{\mathbb{Z}}$ to $\R$. Let $(Y_k)_{k\in\mathbb{Z}}$ be a causal functional of $X$, $i.e.$
$$
\forall n\in\mathbb{Z}, \quad Y_n:=f(\ldots,X_{n-1},X_n).
$$
Then $(Y_n)_{n\geq0}$ is a stationary process with marginal distribution $\nu=\mathcal{L}(Y_0)$. 

The proposition below show that if $(X_n)_{n\in \Z}$ is $\alpha$-mixing ``fast enough" then $H(\theta^*,\cdot)$ ``almost"  lies in $\V_{n^{-\frac 12},2}$ (up to logarithmic factor) as soon as $\E|H(\theta^*,Y_0)|^{2+\delta}<+\infty$ for a $\delta>0$ (so is true for $H(\theta,\cdot)$ since we do not know $\theta^*$ {\em a priori}).

\begin{prop}\label{PropMixingASM} Assume $g\in L^{2+\delta}(\nu)$, $\delta>0$, and that one of the following assumptions holds

\smallskip
\noindent $(a)$  For all $n\in\mathbb{Z}$, $Y_n:=f(\ldots,X_{n-1},X_n)$ where $X$ is a stationary $\alpha$-mixing process satisfying 
\begin{equation}\label{IbragASM}
\sum_{k\ge 1}\frac{\alpha_k^{\frac{\delta}{2(2+\delta)}}}{\sqrt{k}}<+\infty
\end{equation}
and assume the following regularity assumption on $f$
\begin{equation}\label{HypF}
\forall\, 0\leq k'\leq k, \quad \left\|\E\left[g(Y_k)-g(f((0,\ldots,0,X_{k'},\ldots,X_k))\left|\right.\F_0\right]\right\|_2\leq C a_{k-k'}
\end{equation}
for a real constant $C>0$ and a nonnegative real sequence $(a_k)_{k\geq0}$ satisfying
\begin{equation}\label{SommeA}
\sum_{k\ge 1}\frac{a_k}{\sqrt{k}}<+\infty.
\end{equation}
\noindent $(b)$ $Y_n=X_n$, $n\geq0$, and $X$ is a stationary $\alpha$-mixing process satisfying the condition
\begin{equation}\label{IbragimovASM}
\sum_{k\geq 1}\alpha_k^{\frac{\delta}{2+\delta}}<+\infty.
\end{equation}
\begin{equation}\label{GH3ASM}
\mbox{Then }\hskip 3 cm  g \in\V_{\e^{(\eta)}_n,2}, \quad \mbox{with }\e^{(\eta)}_n=\left(\log n\right)^{\frac{3}{2}+\eta}n^{-\frac 12},\mbox{ for every }\eta>0.\hskip 4 cm
\end{equation}
In particular $g$ lies in $\V_{n^{-\beta},2}$  for every  $\beta\in(0,\frac{1}{2})$. 
\end{prop}

\noindent {\bf Remarks.}  Condition~(\ref{IbragimovASM}) is   satisfied when the underlying process $X$ is geometrically $\alpha$-mixing. Slightly refined results could be obtained by calling  upon  Philip and Stout's Law of Iterated Logarithm but the resulting claims would be significantly more technical to state for  little  practical benefit. 

\medskip
An applicationof this result  (based on real data)  is briefly developed in Section~\ref{DkPlASM}. The proof of Proposition\ref{PropMixingASM} relies on the G\`al-Koksma Theorem (see \cite{GalKok} and~\cite{Ben} for a probabilistic version). We state it here in a stationary framework.

\begin{theo}\label{GKThmASM}(G\`al-Koksma's Theorem) Let $\left(\Omega, \F, \P\right)$ be a probability space and let $(Z_n)_{n\geq1}$ be a sequence of random variables belonging to $L^p$, $p\geq1$, satisfying
$$
\E\left|Z_{1}+Z_{2}+\cdots+Z_{N}\right|^p=O(\Psi(N))
$$
where $\frac{\Psi(N)}{N}$, $N\geq1$, is a nondecreasing sequence. Then for every $\eta>0$, 
$$
Z_1(\omega)+Z_2(\omega)+\cdots+Z_N(\omega)=o\left((\Psi(N)(\log(N))^{p+1+\eta})^{\frac{1}{p}}\right)\quad \P(d\omega)\mbox{-}a.s.
$$
\end{theo}

\noindent {\bf Remark.}  The conditions on $X$ and $Z$ come from a result established in \cite{DedMerVol}: by setting $P_0(Z_k):=\E\left[Z_k\left|\right.\F_0\right]-\E\left[Z_k\left|\right.\F_{-1}\right]$, if
\begin{equation}\label{projectionASM}
\sum_{k\in\mathbb{Z}}\left\|P_0(Z_k)\right\|_2<+\infty \quad\mbox{then}\quad\sum_{k\in\mathbb{Z}}\left|\mbox{Cov}(Z_0,Z_k)\right|<+\infty.
\end{equation}
Moreover using \cite{PelUte}, condition (\ref{projectionASM}) is satisfied as soon as 
\begin{equation}\label{CondiASM}
\sum_{k=1}^\infty\frac{1}{\sqrt{k}}\left\|\E\left[Z_k\left|\right.\F_0\right]\right\|_2<+\infty.
\end{equation} 
\paragraph{Proof of Proposition \ref{PropMixingASM}.}  Let $Z_n=g(Y_n)-\int_{{\R}^q} g(y)\nu(dy) $, $n\in\mathbb{Z}$. Without loss of generality we may assume $\int_{{\R}^q} g(y)\nu(dy)=0$.

\smallskip 
\noindent $(a)$ We will rely on  the above Gal-Koksma Theorem (Theorem~\ref{GKThmASM}). First, we evaluate $\E \left|Z_0+\cdots+Z_{n-1}\right|^2$. Setting $S^Z_k=\sum_{j=1}^{k}\E \left[Z_jZ_0\right]$, $k\in \N$,  elementary computations lead to
$$
\E\left|Z_0+\cdots+Z_{n-1}\right|^2=n\E Z_0^2+2\sum_{k=1}^{n-1}\sum_{j=1}^{k}\E\left[Z_{j}Z_0\right] =n\E Z_0^2+2\sum_{k=1}^{n-1}S^Z_{k}=n\left(\E Z_0^2+\frac{2}{n}\sum_{k=1}^{n-1}S^Z_{k}\right).
$$
To establish that $S_n^Z$ converges, we will establish (\ref{CondiASM}).
We set, for $0\leq k'\leq k$, $\F_{k',k}=\sigma(X_{\ell},k'\leq\ell\leq k)$. Then $Z_k=Z_k-\E[Z_k\left.\right|\F_{k',k}]+\E[Z_k\left.\right|\F_{k',k}]$. We derive from the definition of conditional expectation that
$$\left\|\E\left[Z_k-\E[Z_k\left.\right|\F_{k',k}]\left.\right|\F_0\right]\right\|_2\leq \left\|Z_k-\E[Z_k\left.\right|\F_{k',k}]\right\|_2\leq\left\|Z_k-g\circ f(\ldots,0,\ldots,0,X_{k'},\ldots,X_k)\right\|_2$$
so that
\begin{eqnarray*}	\left\|\E\left[Z_k\left|\right.\F_0\right]\right\|_2&\leq&\left\|\E\left[Z_k-\E[Z_k\left.\right|\F_{k',k}]\left|\right.\F_0\right]\right\|_2+\left\|\E\left[\E[Z_k\left.\right|\F_{k',k}]\left|\right.\F_0\right]\right\|_2\\																											&\leq&\left\|\E\left[g\circ f(\ldots,X_0,\ldots,X_{k'},\ldots,X_k)-g\circ f(,\ldots,0,X_{k'},\ldots,X_k)\left|\right.\F_0\right]\right\|_2 \\
                        & & +\left\|\E\left[\E[Z_k\left.\right|\F_{k',k}]\left|\right.\F_0\right]\right\|_2\\
													&\leq& Ca_{k-k'}+\alpha_{k'}^{\frac{1}{r}}\left\|g\circ f(0,\ldots,0,X_{k'},\ldots,X_k)\right\|_p
\end{eqnarray*}
owing to Assumption (\ref{HypF}) on $g\circ f$ and the classical covariance inequality for $\alpha$-mixing process (see \cite{Dou}, Theorem~3(1), p.9) with $\frac{1}{r}+\frac{1}{p}=\frac{1}{2}$, $r,p>2$. As $g\in L^{2+\delta}(\nu)$, $\delta>0$, we may set $p=2+\delta$, $r=\frac{2(2+\delta)}{\delta}$ and $k'=k/2$. Then we obtain
$$\left\|\E\left[Z_k\left|\right.\F_0\right]\right\|_2\leq Ca_{k/2}+\alpha_{k/2}^{\frac{\delta}{2(2+\delta)}}\left\|g\circ f(0,\ldots,0,X_{k/2},\ldots,X_k)\right\|_{2+\delta}.$$
As a consequence, by using (\ref{IbragASM}) and (\ref{SommeA}), we have
$$
\sum_{k=1}^\infty\frac{1}{\sqrt{k}}\left\|\E\left(Z_k\left|\right.\F_0\right)\right\|_2<+\infty,
$$
which implies (owing to   (\ref{projectionASM})) that $S_k^Z$ converges.  Now, by Cesaro's Lemma we have
$$
\E \left| Z_0+\cdots+Z_{n-1}\right|^2=O(n)\quad \mbox{ or equivalently }\quad \left\| \frac 1n\big(Z_0+\cdots+Z_{n-1})\right\|_2=O\Big(n^{-\frac 12}\Big).
$$
Thus, one concludes by by Gal-Koksma's Theorem since,  for every $\eta>0$,
$$
\frac{Z_0+\cdots+Z_{n-1}}{n}=o\left(\frac{\left(\log n\right)^{\frac{3}{2}+\eta}}{\sqrt{n}}\right) \quad  \P\mbox{-}a.s.
$$
\noindent $(b)$ If we assume that $Y_n=X_n$, $n\geq0$, we directly use the covariance inequality for $\alpha$-mixing process
$$
\left|\mbox{Cov} \left(Z_j,Z_0\right)\right|\
\leq8\alpha_{j}^{\frac{1}{r}}\left\|Z_0\right\|_p\left\|Z_0\right\|_q,
$$
where $\frac{1}{r}+\frac{1}{p}+\frac{1}{q}=1$.
By symmetry, we take $p=q>2$ and we get
$$
\left|\E\big(Z_j Z_0\big) \right|\leq8\alpha_{j}^{1-\frac{2}{p}}\left\|Z_0\right\|_p^2. 
$$
As $g\in L^{2+\delta}$, $\delta>0$, we may set $p=2+\delta$ and we obtain $\alpha_{j}^{1-\frac{2}{2+\delta}}=\alpha_{j}^{\frac{\delta}{2+\delta}}$. The condition~(\ref{IbragASM}) can be replaced by the less stringent
Ibragimov's condition~(\ref{IbragimovASM}) to complete the proof. \hfill$\cqfd$

\subsection{Homogeneous Markov chain}

Assume that the innovation process $(Y_n)_{n\geq0}$ is an ${\R}^q$-valued homogeneous Markov chain which transition is $\left(P(y,dz)\right)_{y\in{\R}^q}$ and starting distribution $\mu=\L(Y_0)$. For convenience we will assume that the chain lives on its canonical space $((\R^q)^{\N},\mathcal{B}or(\R^q)^{\otimes\N})$.

\subsubsection{Application of the convergence theorem}

We consider the classical Markov stochastic approximation procedure procedure
\begin{equation}\label{ASMarkovASM}
\th_{n+1}=\th_n-\g_{n+1}K(\th_n,Y_{n+1}), \quad n\geq0,
\end{equation}
where $K:\R^d\times\R^q\to\R^d$ is a Borel function satisfying (\ref{lmrMarkovASM}) below and $\th_0:(\Omega,{\cal A},\P)\to\R^d$ is independent of $(Y_n)_{n\geq0}$. Note that $(Y_n)_{n\ge 0}$ remains is still a Markov chain with respect to   $\F_n=\sigma(\th_0,Y_0,\ldots,Y_n)$, $n\ge 0$.   

Set $H(\th,y):= P\big(K(\theta,.)\big)(y)$ 
and $\Delta M_{n+1}:=K(\th_n,Y_{n+1})-\E\left[K(\th_n,Y_{n+1})\left|\right.\F_n\right]$. 
Then the procedure has the canonical form (\ref{ASpaperASM}) with respect to the filtration $(\F_n)_{n\ge 0}$. 

\medskip
\noindent {\bf Remark.} If we consider that the Markov chain starts from $Y_1$, then $\F_n=\sigma(\th_0,Y_1,\ldots,Y_n)$ and $\E\big[K(\theta_0,Y_1)|\F_0\big]=\E\left[K(\theta,Y_1)\right]_{|\th=\th_0}=\mu P(K(\theta,.))_{|\th=\th_0}$ since $\theta_0$ and $Y_1$ are independent.

\medskip
Let $p\in[1,\infty)$ and set $r=2\vee\frac{p}{p-1}\in[2,+\infty]$. We make the following growth assumption on the function $K$
\begin{equation}\label{lmrMarkovASM}
\forall \th\in{\R}^d, \quad \forall y\in{\R}^q, \quad \left|K(\th,y)\right|\leq C_K\widetilde{\phi}(y)\sqrt{1+L(\th)}
\end{equation}
where $L:\R^d\to\R_+$ satisfies (\ref{LyapunovASM}) and $\sup_{n\geq0}\big\|\widetilde{\phi}(Y_n)\big\|_{r}<+\infty$. 

Then $H$ satisfies (\ref{HHASM}) with $\phi(y)=P\widetilde{\phi}(y)=\E_y\widetilde{\phi}(Y_1)\leq\big\|\widetilde{\phi}\big\|_{L^r(P(y,dz))}<+\infty$ and $\Delta M_{n+1}$ satisfies (\ref{DeltaMASM}) with $\phi(y)=\big\|\widetilde{\phi}\big\|_{L^r(P(y,dz))}$ so that, finally, we may choose $\phi(y)=\big\|\widetilde{\phi}\big\|_{L^r(P(y,dz))}$,  having in mind that $\big\|\phi(Y_n)\big\|_r=\big\|\widetilde{\phi}(Y_{n+1})\big\|_r$.  Now, the  proposition below straightforwardly follows from  Theorem~\ref{Thm1ASM}.
\begin{prop}
Let $p\in[1,\infty)$ and $\th^*\in\R^d$. If $K$ satisfies (\ref{lmrMarkovASM}) and $H$ satisfies the strict pathwise Lyapunov assumption (\ref{lmrASM}), if $(\g_n)_{n\geq1}$ satisfies (\ref{gammaASM}) for a sequence $(\e_n)_{n\geq1}$ satisfying (\ref{HepsASM}) and $H(\th^*,\cdot)\!\in{\cal V}_{\e_n,p}$, then the recursive procedure with Markov innovations defined by (\ref{ASMarkovASM}) converges, $i.e.$
$$\th_n\overset{a.s.}{\underset{n\to+\infty}{\longrightarrow}}\th^*.$$
\end{prop}

\subsubsection{Ergodic framework  description}
We will say that the Markov chain $(Y_n)_{n\geq0}$ (starting from $Y_0\overset{\L}{\sim}\mu$) is $\nu$-ergodic (resp. $\nu$-stable) under $\P_{\mu}$ if for every bounded Borel (resp. continuous) function $f:{\R}^q\rightarrow{\R}$, 
\begin{equation}\label{ErgodicASM}
\P_{\mu}\mbox{-}a.s. \quad \frac{1}{n}\sum_{k=0}^{n-1}f(Y_k)\underset{n\rightarrow\infty}{\longrightarrow}\int_{{\R}^q}fd\nu.
\end{equation}
As soon as the transition $(P(y,dz))_{y\in{\R}^q}$ of $(Y_n)_{n\geq0}$ is Feller, the above $\nu$-stability property implies that $\nu$ is an invariant distribution of the chain, $i.e.$ $\nu P=\nu$. In case of $\nu$-ergodicity the same conclusion holds unconditionally. As a consequence the whole sequence $(Y_n)_{n\geq0}$ is stationary under $\P_{\nu}$.  

\smallskip
Let us focus on the case $\mu=\nu$. If (\ref{ErgodicASM}) holds (with $\mu=\nu$), it is classical background that the whole chain is ergodic under $\P_{\nu}$ (on the canonical space) for the shift operator $\Theta$, $i.e.$ by Birkhoff's theorem, for every functional $F:((\R^q)^{\N},\mathcal{B}or(\R^q)^{\otimes\N})\rightarrow{\R}$, $F\in L^p(\P_{\nu})$,
$$\frac{1}{n}\sum_{k=0}^{n-1}F\circ\Theta^k\underset{n\rightarrow\infty}{\longrightarrow}\E_{\nu}(F) \quad \P_{\nu}\mbox{-}a.s. \mbox{  and in  }L^p(\P_{\nu}),$$
so that by considering $F((y_n)_{n\geq0})=f(y_0)$, $f\in L^p(\nu)$, we finally get that 
$$
{\cal V}_{0^+,p}(\P_{\nu})=L^p(\nu).
$$ 
Note that if the set of invariant distributions for $P$ (convex and) weakly compact and if $\nu$ is extremal in it (so will be $e.g.$  the case  if $\nu$ is unique!) then the chain is ergodic under $\P_{\nu}$ so that the above equality still holds. Furthermore, we know by a straightforward application of G\`al-Koksma Theorem that for any $g\in L^2(\nu)$ for which the related Poisson Equation $g-\nu(g)=\varphi_g-P\varphi_g$ has a solution $\varphi_g\in L^2(\nu)$, then 
\begin{eqnarray*}
\E_{\nu} |g(Y_0)+\cdots+g(Y_{n-1})-n\nu(g)|^2&= &\E_{\nu} |\varphi_g(Y_0)-P\varphi_g(Y_{n-1})+\sum_{1\le k\le n-2}\varphi_g(Y_{k})-P\varphi_g(Y_{k-1})|^2\\
&\le& 6 \nu(\varphi_g^2)+3(n-2)\nu\big((\varphi_g-P\varphi_g)^2\big)=O(n)\\
\mbox{so that  }\hskip 2 cm \bigcap_{\beta\in (0,\frac 12)}\V_{n^{-\beta},2}(\P_{\nu})&\supset &L^2(\nu).\hskip 4 cm
\end{eqnarray*}
Now we will make a connection between the classes $ {\cal V}_{\e_n,p}(\P_\nu)$ and ${\cal V}_{\e_n,p'}(\P_y)$  which will provide examples of non-stationary (Markovian) innovations that can be ``plugged" in stochastic Approximation procedures in the spirit of Theorem~\ref{Thm1ASM}.

\begin{prop} Let $p\in[1,+\infty)$ and let $ p'\in(0, p]$. If $(Y_n)_{n\geq0}$ is $\P_{\nu}$-ergodic and $P(y,dz)=g(y,z)\nu(dz)$, $y\in{\R}^q$ where $g: (\R^q)^2\to  \R_+$ satisfies $\nu(dz)$-$a.e.$, $g(\cdot,z)>0$. Then $\nu$ is the unique invariant distribution of $P$and for every sequence $(\e_n)_{n\geq0}$ satisfying (\ref{HepsASM}),
$$\forall y\in{\R}^q,  \quad g(y,\cdot) \in L^{\frac{p}{p-p'}}(\nu)\Longrightarrow  {\cal V}_{\e_n,p'}(\P_y) \supset {\cal V}_{\e_n,p}(\P_\nu).
$$
\end{prop}
\paragraph{Proof.} It follows form the assumption that any invariant distribution $\nu'$ is equivalent to $\nu$ which implies classically uniqueness.

\smallskip
\noindent$\rhd$ \textit{The $a.s.$ rate.} Let $f\in L^p(\nu)$, $y\in\R^q$ and $\displaystyle 
A_f:=\left\{\omega: \frac 1n \sum_{k=0}^{n-1}f(Y_k(\omega))-\int_{\R^q} fd\nu=O(\e_n)\right\}$.
Since $\liminf_n n \e_n>0$, if $\Theta$ denotes the shift operator on the canonical space of the chain $(Y_n)_{n\geq0}$, $A_f$ clearly satisfies $A_f=\Theta^{-1}(A_f)$ $i.e.$ $\mathds{1}_{A_f}=\mathds{1}_{A_f}\circ\Theta$. Therefore
$$
\P_y(A_f)=\E_y(\mathds{1}_{A_f})=\E_y(\mathds{1}_{A_f}\circ\Theta)=\E_y(\P_{Y_1}(A_f)).
$$
Assume now $f\in\V_{\e_n,p}(\P_{\nu})$. By assumption $\P_{\nu}(A_f)=1$. Let $y\in{\R}^q$. Then 
$$
\P_{\nu}(A_f)=\int_{\R^q}\nu(dz)\P_z(A_f)=1 \quad\mbox{so that}\quad\nu(dz)\mbox{-}a.s. \quad \P_z(A_f)=1.
$$
Now $P(y,dz)\!\ll\!\nu(dz)$ implies $\int_{\R^q}\!P(y,dz)\P_z(A_f)=1$ $i.e.$ $\E_y\left[\P_{Y_1}(A_f)\right]=1$ or equivalently $\P_y(A_f)=1$.

\smallskip
\noindent$\rhd$ \textit{The $L^{p'}$-rate.}  Let $p'\in(0,p]$, let $f \in\V_{\e_n,p}(\P_{\nu})\subset  L^p(\nu)$ and   $\varphi_{n,p'}(y):=\left\|\frac 1n \sum_{k=0}^{n-1}f(Y_k(\omega))-\int_{\R^q} fd\nu\right\|_{L^{p'}(\P_y)}$.   
\begin{equation}\label{HypfLpASM}
\Big\|\frac 1n \sum_{k=0}^{n-1}f(Y_k(\omega))-\int_{\R^q} fd\nu\Big\|_{L^p(\nu)}=O(\e_n) \quad\mbox{so that}\quad\int_{\R^q}\phi^p_n(y)\nu(dy)=O(\e_n^p).
\end{equation}
Assume temporarily that $p'\ge 1$. Consequently, Minkowski inequality implies
\begin{eqnarray*}
	\varphi_{n,p'}(y)&\leq&\frac{\left|f(y)-\int_{\R^q} fd\nu\right|}{n}+\left(1-\frac{1}{n}\right) \left\|\frac{1}{n-1}\sum_{k=1}^{n-1}f(Y_k)-\int_{\R^q} fd\nu\right\|_{L^{p'}(\P_y)}\\
	                     &=&\frac{\left|f(y)-\int_{\R^q} fd\nu\right|}{n}+\left(1-\frac{1}{n}\right)\left(\E_y\left[\E\left(\left|\frac{1}{n-1}\sum_{k=1}^{n-1}f(Y_k)-\int_{\R^q} fd\nu\right|^{p'}\left|\right.\F_1\right)\right]\right)^{\frac{1}{p'}} \\
	                     &=&\frac{\left|f(y)-\int_{\R^q} fd\nu\right|}{n}+\left(1-\frac{1}{n}\right) \left\|\varphi_{n-1,p'}(Y_1)\right\|_{L^{p'}(\P_y)}. 
\end{eqnarray*}
where we used the Markov property in the last equality. Since $P(y,dz)=g(y,z)\nu(dz)$, we derive from H\" older's Inequality (applied to $r=\frac{p}{p'}$ and $s= \frac{p}{p-p'}$) \begin{eqnarray*}
\E_y\varphi_{n-1}(Y_1)^{p'}&=&\int_{\R^q} \varphi_{n-1}(z)^{p'}P(y,dz)=\int_{\R^q} \varphi_{n-1}(z)^{p'}g(y,z)\nu(dz) \\
                        &\leq& \left\|g(y,\cdot)\right\|_{L^{\frac{p}{p-p'}}(\nu)}\left(\int_{\R^q} \varphi_{n-1}(y)^{p}\nu(dy)\right)^{\frac{p'}{p}}\\
                        &\leq&\left\|g(y,\cdot)\right\|_{L^{\frac{p}{p-p'}}(\nu)} \,O(\e_n^{p'}) \quad\mbox{owing to (\ref{HypfLpASM})}.                        
\end{eqnarray*}
$$
\mbox{Finally }\hskip 1.5cm\varphi_n(y)\leq\frac{c}{n}+\left(1-\frac{1}{n}\right)\left\|g(y,\cdot)\right\|_{\frac{p}{p-p'}}\,O(\e_n)= O(\e_n)\quad i.e. \quad f\in\V_{\e_n,p'}(\P_y).\hskip 2cm
$$
The case $p'\in (0,1)$ follows by the usual adjustments (pseudo-Minkowski inequality, etc).$\qquad\hfill\cqfd$ 

\medskip
\noindent {\sc Comments.}  By contrast with the approach of \cite{BMP}, it is not mandatory to solve the Poisson equation related to the pseudo-transition
$$\Pi_{\th_n}(Y_n,dz)=\P\left(Y_{n+1}\in dz\left|\right.\F_n\right)$$
of the algorithm. Indeed, they assume there exists a function $v_{\th}:=v(\th,\cdot)$ solution to
\begin{equation}\label{PoissonASM}
Id-\Pi_{\th}v_{\th}=H(\th,\cdot)-h(\th)
\end{equation}
(Assumption $(H_4)$ in \cite{BMP} p. 220). The target $\th^*$ is then a zero of the mean function $h$ (not canonically defined at this stage in \cite{BMP}). In our setting, $\Pi_{\th}(y,dz)=P(y,dz)$ since the dynamics of $(Y_n)-{n\geq0}$ does not depend upon $\th$, so that Condition (\ref{PoissonASM}) reads
$$v(\th,y)-\int_{\R^q} v(\th,z)P(y,dz)=H(\th,y)-h(\th)$$
where the mean function is naturally defined by
$$h(\th)=\int_{\R^q} H(\th,y)\nu(dy)$$
($\nu$ is the unique invariant probability measure for $P$). Then the family of Poisson equations (indexed by the parameter $\th$) reads
$$v(\th,y)-Pv(\th,y)=H(\th,y)-h(\th).$$
A {\em formal} solution is given by $\displaystyle 
v(\th,y)=\sum_{k\geq0}P^k\left(H(\th,\cdot)-h(\th)\right)(y)$, but  the point is precisely to establish its existence and its properties by using
 the mixing properties of the semi-group $P$ (see~\cite{BMP}).

\section{Applications and numerical examples}\label{cinqASM}

This section is devoted to several examples (mainly inspired by in Finance) of application of convergence theorems in the different frameworks developed in Section \ref{troisASM} and \ref{quatreASM}.

\subsection{Application to implicit correlation search by quasi-stochastic approximation}\label{ICQSA}

Consider a 2-dimensional Black-Scholes model $i.e.$ $X_0^t =e^{rt}$ (riskless asset) and 
$$\forall t\geq0, \quad X^i_t=x^i_0 e^{(r-\frac{\sigma_i^2}{2})t+\sigma_i W^i_t}, \quad x^i_0 > 0, \quad i = 1, 2,$$
for the two risky assets where $\left\langle W^1,W^2\right\rangle_t=\rho t$, $\rho \in [-1, 1]$. Consider a best-of call option characterized by its payoff
$$\left(\max\left(X^1_T ,X^2_T\right)-K\right)_+.$$
We will use a stochastic recursive procedure to solve the inverse problem in $\rho$ 
$$P_{BoC}(x^1_0, x^2_0,K, \sigma_1, \sigma_2, r, \rho, T)=P^{market}_0$$ 
where $P^{market}_0$ is the quoted premium of the option (mark-to-market) with
\begin{eqnarray*}
P_{BoC}(x^1_0, x^2_0,K, \sigma_1, \sigma_2, r, \rho, T)&:=&e^{-rT}\E\left[\left(\max\left(X^1_T,X^2_T\right)-K\right)_+\right] \\
										      &=&e^{-rT}\E\left[\left(\max\left(x^1_0e^{\mu_1T+\sigma_1\sqrt{T}Z^1}, x^2_0
e^{\mu_2T+\sigma_2\sqrt{T}\left(\rho Z^1+\sqrt{1-\rho^2}Z^2\right)}\right)-K\right)_+\right]
\end{eqnarray*}
where $\mu_i=r-\frac{\sigma^2_i}{2}$, $i=1,2$, $Z=(Z^1,Z^2)\stackrel{d}{=}\mathcal{N}(0,I_2)$. We assume from now on that this equation (in $\rho$) has at least one solution, say $\rho^*$. The most convenient way to prevent edge effects due to the fact that $\rho\in[-1, 1]$ is to use a trigonometric parametrization of the correlation by setting $\rho=\cos\th, \th\in{\R}$. This introduces an over-parametrization  since $\th$ and $2\pi-\th$ yield the same solution inside $[0,2\pi]$, but this is not at all a significant problem for practical implementation (a careful examination shows that in fact one equilibrium is repulsive and one is attractive). From now on, for convenience, we will just mention the dependence of the premium function in the variable $\th$, namely
$$\th \longmapsto P(\th):=P_{BoC}(x^1_0, x^2_0,K, \sigma_1, \sigma_2, r, \cos(\th), T).$$
The function $P$ is a $2\pi$-periodic continuous function. Extracting the implicit correlation from the market amounts to solving 
$$P(\th)=P^{market}_0 \quad (\mbox{with} \ \rho=\cos\th).$$
We need the following additional assumption 
$$P^{market}_0 \in (\min_{\th}P,\max_{\th}P)$$
$i.e.$ that $P^{market}_0$ is not an extremal value of $P$. It is natural to set for every $\th\in{\R}$ and every $z=(z^1,z^2)\in{\R}^2$
$$H(\th,z)=e^{-rT}\left(\max\left(x^1_0e^{\mu_1T+\sigma_1\sqrt{T}z^1},x^2_0
e^{\mu_2T+\sigma_2\sqrt{T}(z^1\cos\th+z^2\sin\th)}\right)-K\right)_+- P^{market}_0$$
and to define the recursive procedure 
$$\th_{n+1}=\th_n-\g_{n+1}H(\th_n,Z_{n+1}), \quad n\geq0,\quad \mbox{where} \quad Z_{n+1}\overset{{\cal L}}{\sim}\mathcal{N}(0,I_2),$$
and the gain parameter sequence satisfies (\ref{gammaASM}). For every $z \in {\R}^2$, $\th \longmapsto H(\th,z)$ is continuous and $2\pi$-periodic which implies that the mean function  $h(\th):=\E H(\th,Z_1)=P(\th)-P^{market}_0$ and $\th \longmapsto \E\left[H^2(\th,Z_1)\right]$ are both continuous and $2\pi$-periodic as well (hence bounded).

The main difficulty to apply Theorem \ref{Thm1ASM} is to find out the appropriate Lyapunov function. The quoted value $P^{market}_0$ is not an extremum of the function $P$, hence $\int_0^{2\pi}h^{\pm}(\th)d\th > 0$ where $h^{\pm}:=\max(\pm h, 0)$. We consider $\th_0^*$ any (fixed) solution to the equation $h(\th)=0$ and two real numbers $\beta_{\pm}$ such that
$$0 < \beta_+ <\frac{\int_0^{2\pi}h_+(\th)d\th}{\int_0^{2\pi}h_-(\th)d\th}< \beta_-$$
and we set
$$g(\th):=\left\{
\begin{array}{lcl}
	\mathds{1}_{\{h>0\}}(\th)+\beta_+\mathds{1}_{\{h<0\}}(\th) & \mbox{if} & \th\geq\th_0^* \\
	\mathds{1}_{\{h>0\}}(\th)+\beta_-\mathds{1}_{\{h<0\}}(\th) & \mbox{if} & \th<\th_0^*.\\
\end{array}\right.$$
The function 
$$\th \longmapsto g(\th)h(\th)=h_+-\beta_{\pm}h_-$$
is continuous and ``positively'' $2\pi$-periodic on $[\th_0^*,\infty)$ and ``negatively" $2\pi$-periodic on $(-\infty, \th^*_0]$. Moreover, $gh(\th)=0$ iff $h(\th)=0$ so that $gh(\th^*_0)=gh(\th^*_0-)=0$ which ensures on the way the continuity of $gh$ on $\R$. Furthermore $\int_0^{2\pi}gh(\th)d\th > 0$ and $\int^0_{-2\pi}gh(\th)d\th < 0$
so that, on the one hand,
$$\lim_{\th\rightarrow\pm\infty}\int_0^{\th}gh(u)du = +\infty$$
and, on the other hand, there exists a real constant $C > 0$ such that the function 
$$L(\th)=\int_0^{\th}gh(u)du + C$$
is nonnegative. Its derivative is given by $L'=gh$ so that $L'h=gh^2\geq0$ and $\{L'h=0\}=\{h=0\}$. It remains to prove that $L'$ is Lipschitz continuous. One checks by applying the usual differentiation theorem for functions defined by an integral that, if $\sigma_1\neq\sigma_2$ or $x_1\neq x_2$, then $P$ is differentiable on the whole real line, otherwise it is differentiable only on $\R\setminus 2\pi\Z$, and in both cases
$$
P'(\th)=\sigma_2\sqrt{T}\,\E\left(\mathds{1}_{\{X^2_T>\max(X^1_T,K)\}}X^2_T(\cos(\th)Z^2-\sin(\th)Z^1)\right).
$$
Furthermore, with obvious notations, as soon as $P'(\theta)$ exists, 
$$\left|P'(\theta)\right|\leq\E\left|X^2_T(\cos(\th)Z^2-\sin(\th)Z^1)\right|.$$
The right handside of the inequality defined a $2\pi$-periodic continuous function, hence bounded on the real line. Consequently $\left|P'(\theta)\right|$ is bounded. 
It follows that the $2\pi$-periodic functions $h$ and $h_{\pm}$ are Lipschitz continuous which implies in turn that $L'=gh$ is Lipschitz as well. 

Moreover, one can show that the equation $P(\th)=P^{market}_0$ market has finitely many solutions on every interval of length $2\pi$. One may apply Theorem \ref{Thm1ASM} to derive that $\th_n$ will converge toward a solution $\th^*$ of the equation $P(\th)=P^{market}_0$.  

\medskip
\noindent\textsc{Numerical illustration.} We set the model parameters to the following values
$$x^1_0=x^2_0=100, \ r=0.10, \ \sigma_1=\sigma_2=0.30, \ \rho=-0.50$$
and the payoff parameters
$$
T=1, \ K=100.
$$
The implicit correlation search recursive procedure is implemented with a sequence of some quasi-random normal numbers, namely
$$(\zeta^1_n,\zeta^2_n)=\left(\sqrt{-2\log\left(\xi^1_n\right)}\sin\left(2\pi\xi^2_n\right),\sqrt{-2\log\left(\xi^1_n\right)}\cos\left(2\pi\xi^2_n\right)\right),
$$
where $\xi_n=(\xi^1_n,\xi^2_n), \ n\geq1$, is simply a regular 2-dimensional Halton sequence (see~\cite{Nie} for a definition).

\begin{figure}[!ht]
\centering
\includegraphics[width=15cm]{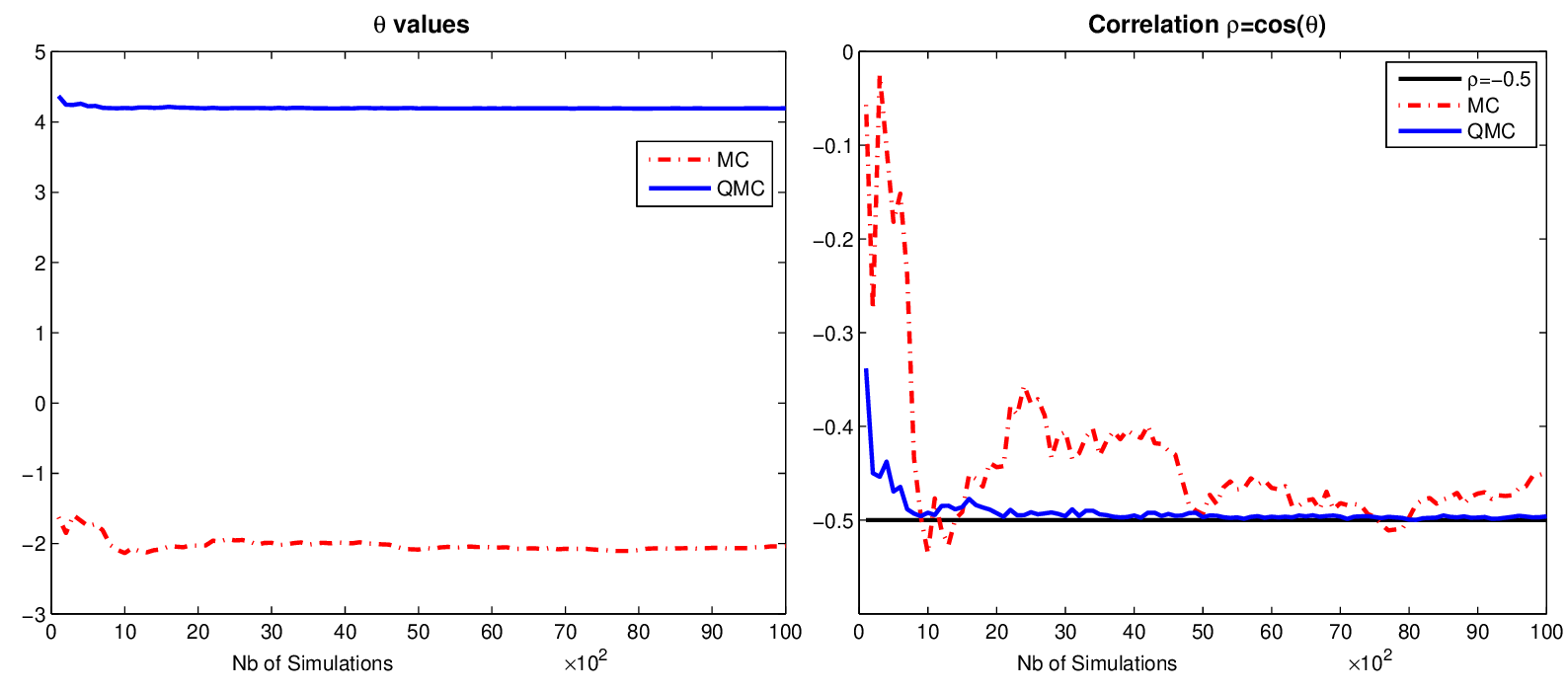}
\vspace{-0.5cm}
\caption{B-S Best-of-Call option. $T = 1$, $r = 0.10$, $\sigma_1 = \sigma_2 = 0.30$, $x^1_0=x^2_0 = 100$, $K = 100$. {\em Left}: convergence of $\th_n$ toward a $\th^*$ (up to $n = 10 000$). {\em Right}: convergence of $\rho_n:= cos(\th_n)$ toward -0.5.}
\label{FigCorrelASM}
\end{figure}
The Black-Scholes reference  price 30.75 is used as a market price so that the target of the stochastic algorithm is $\th^*\in\arccos(-0.5)$. The  parameters of the stochastic approximation procedure are 
$$
\th_0=0, \; n=10^5,\; \g_n=\frac{8}{n},\; n\ge 1.
$$

The choice of $\th_0$ is ``blind'' on purpose (see~Figure~\ref{FigCorrelASM}).

\subsection{Recursive omputation of the VaR and the CVaR}\label{VaRCVaR}

Another example of application is the recursive computation of financial risk measure which are the best known and the most common: the Value-at-Risk (VaR) and the Conditional Value-at-Risk (CVaR). This risk measures evaluate the extreme losses of a portfolio potentially faced by traders. The recursive computation of the VaR and the CVaR was introduced in \cite{BarFriPag}, based on the formulation as an optimization problem (see \cite{RocUry}) and on an unconstrained importance sampling procedure developed in \cite{LemPag}. These  variance reduction aspects are not investigated here.

\subsubsection{Definitions and formulation}
Let $Y:(\Omega,A,\P)\rightarrow{\R}$ be a random variable representative of a loss ($Y\geq0$ is a loss equal to $Y$).
\begin{defi}  The Value at Risk (at confidence level $\alpha\in(0,1)$, $\a\approx 1$) of a given portfolio is the (lowest) $\alpha$-quantile of the distribution $Y$ $i.e.$
$$VaR_{\alpha}(Y):=\inf\left\{\th\,|\,\P(Y\leq\th)\geq\alpha\right\}.$$
\end{defi}
As soon as the distribution function of $Y$ has no atom, the value at risk satisfies $P(Y\leq VaR_{\alpha}(Y))=\alpha$ and if the distribution function $F_Y$ of $Y$ is also increasing (strictly) then, it is the unique solution. As this risk measure is not consistent (see \cite{FolSch}), another consistent risk measure is provided by the Conditional value at Risk when $Y\in L^1(\P)$ with a continuous distribution (no atom).
\begin{defi}
Let $Y\in L^1(\P)$ with an atomless distribution. The Conditional value at Risk (at level $\alpha$) is the conditional expectation of the portfolio loss $Y$ beyond $VaR_{\alpha}(Y)$, $i.e.$
$$CVaR_{\alpha}(Y):=\E\left[Y|Y\geq VaR_{\alpha}(Y)\right].$$
\end{defi}
The following formulation of the $VaR_{\alpha}(Y)$ and $CVaR_{\alpha}(Y)$ as solutions to an optimization problem is due to Rockafellar and Uryasev in \cite{RocUry}.
\begin{prop}
(Rockafellar and Uryasev) Let $Y\in L^1(\P)$ with an atomless distribution. The function $V:\th\mapsto\th+\frac{1}{1-\alpha}\E\left(Y-\th\right)_+$ is convex, and
$$CVaR_{\alpha}(Y)=\min_{\th}\left(\th+\frac{1}{1-\alpha}\E\left(Y-\th\right)_+\right) \,\, \mbox{with} \,\, VaR_{\alpha}(Y)=\inf \arg\min_{\th}\left(\th+\frac{1}{1-\alpha}\E(Y-\th)_+\right).$$
\end{prop}

\subsubsection{Stochastic gradient for the computation of both  \texorpdfstring{$VaR_{\alpha}(Y)$}{VaR(Y)} and  \texorpdfstring{$CVaR_{\alpha}(Y)$}{CVaR(Y)}}

$\rhd$  {\em Computation of the \texorpdfstring{$VaR_{\alpha}(Y)$}{VaR(Y)}.} What precedes  suggests to implement a stochastic gradient descent derived from the above {\em convex} objective function $V(\th)=\th+\frac{1}{1-\alpha}\E(Y-\th)_+$. Assume that $Y\in L^1(\P)$ with a continuous increasing distribution function $F_Y$ (for the sake of simplicity, see~\cite{BarFriPag} for a slightly more general framework). 
Let $\nu={\cal L}(Y)$. We check that
$$
\lim_{\th\rightarrow+\infty}\frac{V(\th)}{\th}=1 \quad \mbox{and} \quad \lim_{\th\rightarrow+\infty}\frac{V(-\th)}{\th}=\frac{\alpha}{1-\alpha} \quad \mbox{hence} \quad \lim_{\th\rightarrow\pm\infty}V(\th)=+\infty.
$$
so that $\{VaR_{\alpha}(Y)\}=\mathop{{\rm argmin}}_{\R}V$. We check that $V'(\th)=\E[H(\th,Y)]$ where 
$$
H(\th,y):=1-\frac{1}{1-\alpha}\mathds{1}_{\{y\geq\th\}}.
$$
Note that $H$ is uniformly bounded by $1\vee\frac{\alpha}{1-\alpha}$. This leads to  devise the stochastic gradient descent
$$
\th_{n+1}=\th_n-\g_{n+1}H(\th_n,Y_{n}), \quad n\geq0, \quad \th_0\in L^1(\P).
$$
whose unique  target is $\theta^*=VaR_{\alpha}(Y)$. It is clear that, for every $y\in\R$, $\theta\mapsto H(\theta,y)$ is nondecreasing so that $L(\theta)=\frac 12(\theta-\theta^*)^2$ is a good candidate as a Lyapunov function. In fact it is even a strict pathwise Lyapunov function in the sense of~(\ref{lmrASM}) by setting for every $\delta>0$, $\Psi_{\delta}(\theta)= \delta \mathds{1}_{|\theta-\theta^*|>\delta}$ and $\chi_{\delta}(y) =\mathds{1}_{|y-\theta^*|\le \delta}$. 

\smallskip As soon as $(Y_n)_{\ge0}$ is $\nu$-averaging, there exists a sequence $(\e_n)_{n\ge 1}$ such that, for every $\theta\in \R$, $\mathds{1}_{\{y\ge \theta\}}\in {\cal V}_{\e_n,2}$ since the empirical distribution measure $a.s.$ (and subsequently in $L^2$)  converges uniformly toward $F_Y$. 
Finally, as soon as the step sequence $(\gamma_n)_{n\ge 1}$ is admissible for $(\e_n)_{n\ge 1}$, Theorem \ref{Thm1ASM} implies that
$$\th_n\overset{a.s.}{\underset{n\rightarrow\infty}{\longrightarrow}}\th^*=VaR_{\alpha}(Y).$$
In practice $\gamma_n=c/n$, $c>0$, is always admissible given the rate of convergence of the empirical measure in usual applications. Of course, when the $Y_n$ are i.i.d., standard martingale arguments ``\`a la Robbins-Monro" make things straightforward under less stringent assumptions on the step sequence.

\smallskip
\noindent $\rhd$ {\em  Computation of the \texorpdfstring{$CVaR_{\alpha}(Y)$}{CVaR(Y)}}. The idea to compute the $CVaR_{\alpha}(Y)$ is to devise a companion procedure of the above stochastic gradient by setting, $\zeta_0=0$  and for every $n\ge 0$, 
$$
\zeta_{n+1}=\zeta_n-\frac{1}{n+1}\left(\zeta_n-v(\th_n,Y_{n})\right) \quad \mbox{ with }\quad v(\th,y):=\th+\frac{(y-\th)_+}{1-\alpha}.
$$
One checks that for every $n\ge 0$, 
\[
\zeta_n = \frac 1n \sum_{k=0}^{n-1} v(\theta_k,Y_{k})=  \frac 1n \sum_{k=0}^{n-1} v(\theta^*,Y_{k})+ \frac 1n \sum_{k=0}^{n-1} v(\theta_k,Y_{k})- v(\theta^*,Y_{k})
\]
Using that $v$ is Lipschitz continuous in $\theta$ uniformly in $y$, we derive that the second term in the right hand side of the above equality goes to $0$  $a.s.$ as $\theta_n\to 0$ $a.s.$  

As concerns the first term,  still in right hand side,  first note that $v(\theta^*,y)$ has a linear growth in $y$ so  it will $a.s.$ go to $\E\, v(\theta^*,Y)= V(\theta^*)= CVaR_{\alpha}(Y)$    as soon as, $e.g.$, $\sup_{n\ge 1}\frac1n\sum_{k=0}^{n-1}|Y_k|^{1+\eta}<+\infty$ $a.s.$ for an $\eta>0$ by combining standard uniform integrability arguments (with respect to the empirical measure) and the $\nu$-stability of $(Y_n)_{n\ge 0}$.
In practice one must keep in mind that an  adaptive importance sampling procedure like that detailed in \cite{BarFriPag} should be added. For a $QMC$ implementation of the procedure, see~\cite{Fri}.

\subsection{Long term investment evaluation (and inhomogeneous Markov innovations)}\label{LTIIMI}

In this example we deal with  averaging {\em inhomogeneous} Markov innovations, namely the Euler scheme with decreasing step of a Brownian diffusion. To describe the functional class   $\V_{\e_n,p}$, we rely on an approach  developed in \cite{LamPag} and \cite{Lem} to compute the invariant measure of a diffusion.

\subsubsection{Computation of the invariant distribution of a diffusion}

We consider a stochastic recursive algorithm for the computation of the invariant distribution $\nu$ introduced in \cite{LamPag} of a Brownian diffusion process 
\begin{equation}\label{diffusionASM}
dY_t=b(Y_t)dt+\sigma(Y_t)dW_t
\end{equation}
where $b:{\R}^q\rightarrow{\R}^q$  and  $\sigma: \R^q \to \mathcal{M}_{q,\ell}(\R)$ (matrices with $q$ rows and $\ell$ columns) are Lipschitz continuous, and $W$ is a $\ell$-dimensional Brownian motion. We denote by $\mathcal{A}$ its infinitesimal generator and by $(P_t)_{t\geq0}$ its transition semi-group.

First, we introduce the Euler discretization of (\ref{diffusionASM}) with a step $\g_n$ vanishing to 0, $i.e.$
\begin{equation}\label{EulerASM}
\forall n\in\mathbb{N}, \quad \bar{Y}_{n+1}=\bar{Y}_n+\bar{\g}_{n+1}b(\bar{Y}_n)+\sqrt{\bar{\g}_{n+1}}\sigma(\bar{Y}_n)U_{n+1},
\end{equation}
where $\bar{Y}_0\in L^0_{{\R}^q}(\Omega,\F,\P)$ and $(U_n)_{n\geq1}$ is ${\R}^{\ell}$-valued normalized white noise defined on a probability space $(\Omega,\F,\P)$, independent of $\bar{Y}_0$. The step sequence $\bar{\g}:=(\bar{\g}_n)_{n\geq1}$ satisfies the conditions
\begin{equation}\label{GamASM}
\forall n\geq1, \quad \bar{\g}_{n}\geq0, \quad \lim_{n\rightarrow\infty}\bar{\g}_{n}=0 \quad \mbox{and} \quad \bar{\Gamma}_{n}:=\sum_{k=1}^n\bar{\g}_k\underset{n\rightarrow\infty}{\longrightarrow}+\infty.
\end{equation}
\noindent 
For every $n\geq1$ and every $\omega\in\Omega$, set
\vskip-0.5cm
\begin{equation}\label{NuASM}	
	\nu_n(\omega,dy):=\frac{1}{n}\sum_{k=0}^{n-1}\delta_{\bar{Y}_k(\omega)}.
\end{equation}
We will use $\nu_n(\omega,f)$ which can be compute recursively to approximate $\nu(f)$.
\begin{defi}(Strong condition of stability) A diffusion with generator $\mathcal{A}$ satisfies a strong stability condition of type $(V,\alpha)$ if there exists a (so-called Lyapunov) function $V\in\mathcal{C}^2({\R}^q,[1,+\infty[)$ such that $\lim_{|y|\rightarrow+\infty}V(y)=+\infty$ and $\exists\,\alpha>0$, $\exists\,\beta\geq0$ such that $\mathcal{A}V\leq -\alpha V+\beta$.
\end{defi}

\noindent {\bf Remark.}  If the $(V,\alpha)$-strong stability condition holds then (\ref{diffusionASM}) admits a strong solution starting from any $y\in\R^d$ and admits at least one invariant distribution $\nu$ ($i.e.$ $\nu P_t=\nu$, $t\geq0$).
\begin{defi} $(a)$ A couple $(\bar{\gamma},\eta)$ is an {\rm averaging step-weight system} if the sequences $(\bar{\gamma}_n)_{n\geq1}$ and $(\eta_n)_{n\geq1}$ are nonnegative, general terms of a non-converging series and such that
$$
\lim_{n}\bar{\gamma}_n=0, \quad \sum_{n\geq1}\frac{1}{H_n}\left(\Delta \frac{\eta_n}{\bar{\g}_n}\right)_+<+\infty \quad \mbox{and} \quad  \sum_{n\geq1}\left(\frac{\eta_n}{H_n\sqrt{\bar{\g}_n}}\right)^2<+\infty,
$$
where $H_n=\sum_{k=1}^n\eta_k$.

\noindent $(b)$ In particular, if $\eta_n\equiv1$, then $(\bar{\gamma}, 1)$ is an {\rm averaging step-weight system} if
$$\lim_{n}\bar{\gamma}_n=0, \quad \sum_{n\geq1}\frac{1}{n}\left(\frac{1}{\bar{\g}_n}-\frac{1}{\bar{\g}_{n+1}}\right)<+\infty \quad \mbox{and} \quad  \sum_{n\geq1}\frac{1}{n^2\bar{\g}_n}<+\infty.$$
\end{defi}

The terminology ``averaging'' refers here to the fact that if ${\cal A}$ is $(V,\alpha)$-stable (and the invariant distribution $\nu$ is unique  for the sake of simplicity) then, as soon as $(\bar{\gamma}, \eta)$ is averaging (see $e.g.$ \cite{LamPag}, \cite{Lem} or \cite{Lem2}), then 
\vskip-0.5cm
$$\sup_{n\geq0} \E V(\bar Y_n)<+\infty\quad \mbox{ and } \quad \P(d\omega)\mbox{-}a.s. \quad \nu^{\eta}_n(\omega,dy):=\frac{1}{H_n}\sum_{k=0}^{n-1} \eta_k\delta_{\bar{Y}_k(\omega)}\overset{(\R^d)}{\Longrightarrow}\nu.$$

\noindent {\bf Example.} If $\bar{\g}_n=\frac{\bar{\g}_0}{n^r}$, $0<r<1$, and $\eta_n\equiv1$, then $(\bar{\gamma},1)$ is averaging.  

\smallskip
We assume that the diffusion $(Y_t)_{t\geq0}$ satisfies a strong condition of stability of type $(V,\alpha)$ with $V$ sub-quadratic and that the invariant measure $\nu$ is unique. Besides the coefficients $b$ and $\sigma$ satisfy $|b|^2+\mbox{Tr}(\sigma\sigma^t)=O(V)$. Then the Euler scheme with decreasing step $(\bar{Y}_n)_{n\geq0}$ defined by (\ref{EulerASM}) satisfies a strong condition of stability of type $(W,n_0)$ where $W$ is a function depending upon $V$ and the moments of $U_1$, namely
$$
\forall n\geq n_0, \quad \E\left[W(\bar{Y}_{n+1})\left.\right|\sigma\left(\bar{Y}_k, 0\leq k\leq n\right)\right]\leq\left(1-\alpha\,\bar{\gamma}_{n+1}\right)W(\bar{Y}_n)+\beta.
$$

$\rhd$ If $U_1$ is sub-normal (typically if $U_1$ {\em is} normal), and $\Tr(\sigma\sigma^*)\leq C_{\sigma} V^{1-\zeta}$, we may choose $W=\exp(\lambda V^{\zeta})$ with $\lambda$ small enough (see \cite{Lem} Proposition III.2 p. 36).  

\smallskip
$\rhd$ If $U_1$ has a moment of order $2(p+1)$, $p\geq2$, then $W=V^{p+1}$ ({\em idem}). 

\medskip
Assume that the function $f:\R^q\rightarrow\R$ admits a regular enough solution $\phi$ to the Poisson equation
\begin{equation}\label{poisson2ASM}
\mathcal{A}\phi=-(f-\nu(f)),
\end{equation}
 $i.e.$ belonging to the set
$$
\mathcal{E}_{p,W}:=\left\{\phi\in\mathcal{C}^p({\R}^q,{\R}), \, \forall j\in\{0,\ldots,p\}, \, \forall y\in{\R}^q, \, \left|D^j\phi(y)\right|^2=o\left(\frac{W(y)}{V^j(y)}\right)\right\}
$$
and satisfying $D^p\phi$ Lipschitz. For such functions $\phi$, let us to define the functions $D_q$, $3\leq q\leq p$, by
$$
\forall y\in\R^q, \quad D_q(y)=\sum_{j\geq q/2}^q\frac{C_j^{q-j}}{j!}D^j\phi(y)\cdot\left(\left(b(y)\right)^{\otimes(q-j)},\E\left[\left(\sigma(y)U_1\right)^{\otimes(2j-q)}\right]\right).
$$
They will appear in the development of the error of order $p$, $p\geq3$.

\begin{theo}\label{L2CondASM}
Let $p\geq2$ such that $U_1\in L^{2(p+1)}$ and $\phi\in\mathcal{E}_{p,W}$ solution to Poisson equation (\ref{poisson2ASM}) such that $D^p\phi$ is Lipschitz. Define $q^*$ by
$$q^*=\min_{q\in\{3,\ldots,p\}}\{D_q\neq0\}\wedge(p+1).$$
(Note that if $U_1\overset{{\cal L}}{\sim}{\cal N}(0,I_q)$ then $q^*=4$). Let $\bar{\Gamma}_n^{(\beta)}=\sum_{k=1}^n\bar{\g}_k^{\beta}$, $\beta\in\R$. Assume that the couples $(\bar{\gamma},1)$ and $(\bar{\gamma},\frac{1}{\bar{\g}})$ are averaging and that $(\bar{\g}_n)_{n\geq1}$ is non-increasing. If $ q^*\leq p$ and 
$$
\frac{\bar{\Gamma}_n^{(q^*/2-1)}}{\sqrt{\bar{\Gamma}_n^{(-1)}}}\underset{n\rightarrow+\infty}{\longrightarrow}\xi\in]0,+\infty], \quad\quad \Big(\Big(\bar{\Gamma}_n^{(q^*/2-1)}\bar{\g}_n\Big)^{-1}\Big)_{n\geq1} \mbox{ is non-increasing},$$
$$
\sum_{n\geq1}\frac{1}{\bar{\Gamma}_n^{(q^*/2-1)}}\left|\Delta \frac{1}{\bar{\g}_n}\right|<+\infty \quad \mbox{and} \quad \sum_{n\geq1}\frac{1}{\bar{\g}_n\left(\bar{\Gamma}_n^{(q^*/2-1)}\right)^2}<+\infty,
$$
then 
$$
f\in\V_{\e_n,2}\quad \mbox{with} \quad \e_n=\frac{\bar{\Gamma}_n^{(q^*/2-1)}}{n}\underset{n\rightarrow+\infty}{\longrightarrow}0.
$$
\end{theo}

\begin{cor}\label{epsEulerInvASM}
If $\bar{\g}_n=\frac{\bar{\g}_0}{n^r}$, $\bar{\g}_0>0$, the above theorem holds true when $0<r\leq\frac{1}{q^*-1}$ and $\e_n=n^{-r\left(q^*/2-1\right)}$. In particular, for a  Gaussian Euler scheme, $\e_n =n^{-r}$. 
\end{cor}

\noindent {\bf Sketch of proof of Theorem \ref{L2CondASM}.} Proposition V.4 in \cite{Lem} gives
$$
\Big\|\frac{n}{\bar{\Gamma}_n^{(q^*/2-1)}}\Big(\frac{1}{n}\sum_{k=0}^{n-1}f(\bar{Y}_k)-\nu(f)\Big)\Big\|_2=\left\|M_n+S_n\right\|_2+o(1)
$$
$$
\mbox{where }\; M_n=\frac{1}{\bar{\Gamma}_n^{(q^*/2-1)}}\sum_{k=1}^n\sqrt{\frac{1}{\bar{\g}_k}}\left\langle \nabla\phi(\bar{Y}_{k-1})\left.\right|\sigma(\bar{Y}_{k-1})U_k\right\rangle\;\mbox{ and }\; S _n=\sum_{q=q^*}^p\frac{1}{\bar{\Gamma}_n^{(q^*/2-1)}}\sum_{k=1}^n\bar{\g}_k^{\frac{q}{2}-1}D_q(\bar{Y}_{k-1}).
$$
Using that $\phi\in\mathcal{E}_{p,W}$, $i.e.$ that for every $q$, $\left|D_q\right|^2=o(W)$, and that $\sup_n\E\,W(\bar{Y}_n)<+\infty$ (according to the stability condition of the Euler scheme), we get (see the  remark after Proposition V.1 p.62 in \cite{Lem}), 
$$
\left\|M_n\right\|_2\leq\frac{\sqrt{\bar{\Gamma}_n^{(-1)}}}{\bar{\Gamma}_n^{(q^*/2-1)}}\sup_{0\leq k\leq n-1}\left\|\sigma^*\nabla\phi(\bar{Y}_k)\right\|_2<+\infty
$$
since $\sqrt{\bar{\Gamma}_n^{(-1)}}/\bar{\Gamma}_n^{(q^*/2-1)}\underset{n\rightarrow+\infty}{\longrightarrow}\xi^{-1}\in]0,+\infty[$ and $\displaystyle  \|\sup_n S_n\|_2<+\infty$. \hfill$\cqfd$

\subsubsection{Application to the minimization of a potential}
The aim is to minimize  a convex potential $V:\R^q\to\R$ having a  minimum ($e.g.$ because $\lim_{|\theta|\to+\infty}V(\theta)=+\infty$) assumed to be unique. We also assume that $V$ has a representation as the expectation with respect to the invariant distribution $\nu$ of an ergodic diffusion, say $Y$ defined above. Typically $V$ appears as the long run limit (under appropriate assumptions) of a functional through Birkhoff's Theorem:
\[
V(\theta)= \lim_{t\to+\infty} \frac1T\int_0^T v(\theta,Y_{t+s})\,ds= \E_{\nu}\Big(\frac 1T \int_0^Tv(\theta,Y_{s}\,ds\Big)= \int_{\R^q}v(\theta,y)\nu(dy).
\]
We make the following assumptions
\begin{itemize}
	\item[$(i)$] Integrability: $\forall y\in{\R}^q, \;\th\mapsto v(\th,y)$ is convex.
	\item[$(ii)$] Pathwise convexity: $\forall \th\in{\R}^d, \;v(\th,\cdot)\in L^1(\nu)$.
	\item[$(iii)$]  Differentiability: $\forall \th\in{\R}^d, \; \nabla_{\!\th}v(\th,y)$ exists.
	\item[$(iv)$]  Uniform integrability:  $\forall \th\!\in{\R}^d, \,\Big(\frac{| v(\th,y)-v(\th',y)|}{|\th-\th'|}\Big)_{\theta'\in[\th-\eta_\th,\th+\eta_\th]\setminus\{\th\}}$, $\eta_\th\!>\!0$, is uniformly integrable.
\end{itemize}
Then (using uniqueness of $\th^*$),
$$
\theta^*={\rm argmin}_{\th\in{\R}^d}\int_{{\R}^q}  v(\th,y)\nu(dy)\; \mbox{ iff }\;\int_{{\R}^q} \nabla_{\th}v(\th^*,y)\nu(dy)=0,
$$
At this stage the idea is to devise a stochastic gradient (gradient based recursive zero search) using the Gaussian Euler scheme $(\bar Y_n)_{n\ge 0}$ with decreasing step $\bar{\g}_n =\bar{\g}_0n^{-\frac 13}$, $\bar{\g}_0>0$, of $Y$ as an $\nu$-averaging innovation process with rate $\e_n = \bar{\Gamma}_n/n\to 0$:
\[
\forall n\geq0, \quad \th_{n+1}=\th_{n}-\g_{n+1}\nabla_{\!\th}v(\th_n,\bar{Y}_n),\; \th_0\!\in \R^q.
\]
Let $p\!\in[1,\infty)$ such that $\nabla_{\!\th} v$ satisfies the growth assumption~(\ref{HHASM}) with $L(\th)=|\th-\th^*|^2$, $\nabla_{\!\th}v(\th^*,\cdot)\!\in\V_{\e_n,p}$ and $(\g_{n})_{n\ge 1}$ is admissible for $\e_n$ given by Corollary~\ref{epsEulerInvASM}, then Theorem~\ref{Thm1ASM} implies that $\theta_n \to\theta^*$~$a.s.$

\medskip
\noindent {\sc Toy numerical example.} We consider a long-term investment project (see the example in~\cite{LokZer}) which yields payoff at a rate that depends on the installed capacity level and on the value of an underlying state process modeled with an ergodic diffusion. The process $Y$   represents an economic indicator such as the asset demand or its discounted price. Our aim is to determine the capacity expansion strategy that {\em maximizes}  the long-term average payoff resulting from the project operation. So it is an ergodic control problem in a microeconomic framework. In~\cite{LokZer} is shown that this dynamical optimization problem is  equivalent (see above) to a static optimization problem involving the stationary distribution $\nu$ of $Y$ and the (concave) running payoff function $C$, namely, still following~\cite{LokZer}, 
$$
\forall \,\th\!\in{\R}_+, \;\forall \,y\!\in\R_+, \quad C(\th,y)=y^{\alpha}\th^{\beta}-c\,\th \quad \mbox{where $\;\alpha,\beta\in(0,1)\;$ and $\;c\in(0,\infty)$.}
$$
The term $y^{\alpha}\th^{\beta}$ can be identified to the so-called {\em Cobb-Douglas production function}, while the term $c\,\th$  measures  the cost of capital use.  Our task is to {\em minimize} $\displaystyle \int_{\R^q}(-1)\big(y^\a\th^\b-c\,\th\big)\nu(dy)$ (so that of course  $\th^*= \big(\frac{\b \E Y_1^\a}{c}\big)^{\frac{1}{1-\b}}$!). Since $\nabla_{\!\th} \,C(\th,y)$ is singular at $\th=0$, we will introduce the increasing convex with linear growth change of variable $\th=\big(\tilde\th+(\tilde\th^2+1)^{1/2}\big)^{\rho(\tilde\th)}$, $\rho(\tilde\th)=\frac{1}{1-\beta}\mathds{1}_{\tilde\th<0}+\mathds{1}_{\tilde\th\geq0}$, from $\R$ onto $(0,\infty)$ and we  consider
 \[
  \nabla_{\!\th}v(\widetilde \th,y)= -\nabla_{\!\th} C\Big(\big(\tilde\th+(\tilde\th^2+1)^{1/2}\big)^{\rho(\tilde\th)},y\Big) , \;\tilde \th\!\in \R, \;y\!\in\R_+.
 \]
Still following~\cite{LokZer}, the dynamics of the underlying state process $Y$ is modeled by the one-dimensional CIR diffusion (whose diffusion coefficient is unfortunately not Lipschitz), namely
\begin{equation}
dY_t=\kappa\left(\vartheta-Y_t\right)dt+\sigma\sqrt{\left|Y_t\right|} dW_t, \quad Y_0>0,
\end{equation}
where $\kappa,\vartheta,\sigma>0$ are constants satisfying $2\kappa\vartheta>\sigma^2$ so that $(Y_t)_{t\geq0}$ is $(0,\infty)$-valued. 

The resulting stochastic gradient procedure with step $(\widetilde \g_n)_{n\ge 1}$ reads
$$
\forall n\geq0, \quad \widetilde\th_{n+1}=\widetilde\th_{n}-\g_{n+1} \nabla_{\!\th}v(\widetilde \th_n,\bar Y_{n}),\quad \widetilde \th_0 \!\in \R,
$$
where $(\g_n)_{n\ge1}$ is admissible with respect to $\e_n =  \bar{\Gamma}_n/n$  and $(\bar Y_n)_{n\ge0}$ the Euler scheme with step $\bar{\g}_n=\bar{\g}_0n^{-1/3}$ ($L(\widetilde\th)=|\widetilde \th-\widetilde \th^*|^2$ is still a pathwise Lyapunov function). One checks that $\nabla_{\!\th}v$ satisfies~(\ref{HHASM}) with $\phi(y)\equiv c_{\a,\b} y^\a$, $c_{\a,\b}>0$ and $p=2$ since $\sup_n \E\, \bar Y_n^2<+\infty$ and $\a\!\in(0,1)$.

\begin{figure}[!ht]
\vspace{-0.25cm}
\centering
\includegraphics[width=10cm]{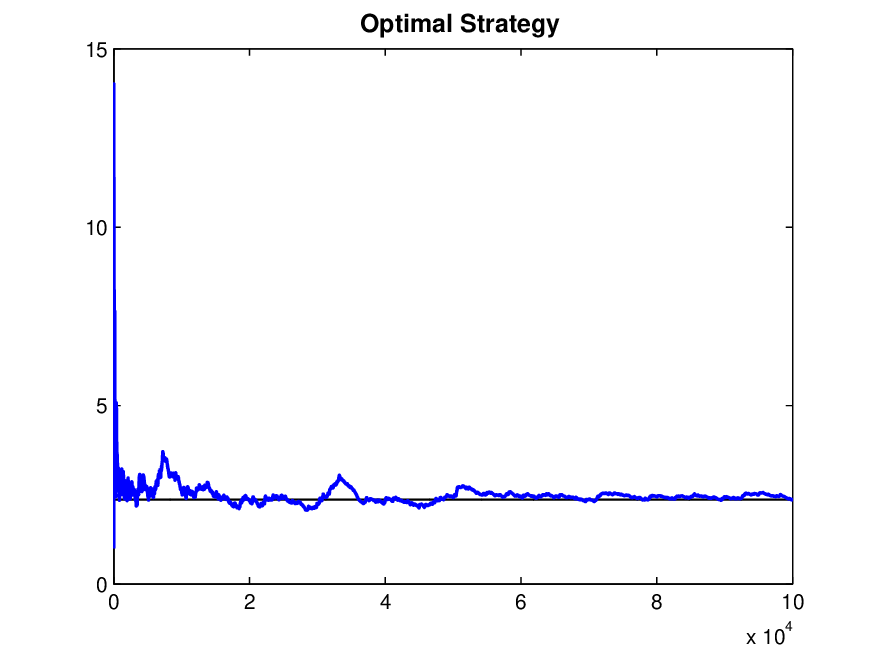}
\vspace{-0.5cm}
\caption{Convergence towards the optimal capacity level of the investment project: $\kappa=1$, $\vartheta=1$, $\sigma=1.5$, $\alpha=0.8$, $\beta=0.7$, $c=0.5$, $n=10^5$, $\widetilde \g_n=\frac{5}{n}$.}
\label{FigInvariantDistribASM}
\end{figure}

\smallskip The invariant distribution of $Y$ is a Gamma law which density is given by
$$
\nu(dy)=\frac{1}{\Gamma\left(\frac{2\kappa\vartheta}{\sigma^2}\right)}y^{\frac{2\kappa\vartheta}{\sigma^2}-1}\exp\left(\frac{2\kappa}{\sigma^2}\left[\vartheta\log\left(\frac{2\kappa}{\sigma^2}-y\right)\right]\right)\mathds{1}_{\{y>0\}},
$$
where $\Gamma$ is the gamma function. Thus we can compute the previous integral, namely
$$
\int_{{\R_+}}y^{\alpha}\nu(dy)=\frac{\Gamma\left(\frac{2\kappa\vartheta}{\sigma^2}+\alpha\right)}{\Gamma\left(\frac{2\kappa\vartheta}{\sigma^2}\right)}\left(\frac{\sigma^2}{2\kappa}\right)^{\alpha}<+\infty,
$$
so we have in fact  a closed form for $\th^*$ given by $\theta^*=\left(\displaystyle\frac{\beta\Gamma\left(\frac{2\kappa\vartheta}{\sigma^2}+\alpha\right)}{c\Gamma\left(\frac{2\kappa\vartheta}{\sigma^2}\right)}\left(\frac{\sigma^2}{2\kappa}\right)^{\alpha}\right)^{\frac{1}{1-\beta}}$. Figure~\ref{FigInvariantDistribASM} illustrates the convergence of the algorithm (the parameters are specified in the caption).

If one considers a basket of assets modeled by a Wishart process (see \cite{Bru} and \cite{GouJasSuf}), a similar long-term ergodic control process can be devised. Closed forms are no longer available for the static optimization problem. However, our numerical approach can be extended straightforwardly (provided one has at hand an efficient method of simulation for Wishart process, like that proposed in \cite{GouJasSuf}).

\subsection{The ergodic two-armed bandit} 
An application of the multiplicative setting  is the so-called two-armed bandit algorithm introduced in mathematical psychology, learning automata (see \cite{Nor,NarSha}) and more recently asset allocation (\cite{LamPagTar}). Criteria on $a.s.$ convergence under pure  i.i.d. assumptions were obtained in \cite{LamPagTar, LamPag2} and under ergodic assumptions in \cite{TarVan}. A penalized version of this algorithm is also studied in~\cite{LamPag3}. 

This algorithm is defined as follows: at each step $n\geq0$, one plays arm $A$ (resp. arm $B$) at random with probability $\th_n$ (resp. $1-\th_n$), where $\th_0=\th\in(0,1)$ and $\th_n$ is updated according the following ``{\em rewarding}" rule:  for every  $n\geq0$,
\begin{equation}\label{BanditASM}
\th_{n+1}=\th_n+\g_{n+1}\Big((1-\th_n)\mathds{1}_{\{U_{n+1}\leq\th_n\}\cap A_{n+1}}-\th_n\mathds{1}_{\{U_{n+1}>\th_n\}\cap B_{n+1}}\Big)
\end{equation}
where $(U_n)_{n\geq1}$ is an i.i.d. sequence of uniform random variables, independent of $(A_n)_{n\geq1}$ and $(B_n)_{n\geq1}$ which  are two  sequences of (possibly dependent) events evaluating the performances of the arms $A$ and $B$ respectively ($A_n$ is the event ``$A$'s performance is satisfactory ay time $n$" idem for $B_n$ and $B$). 

This stochastic procedure can be rewritten in a canonical form as follows
\begin{equation}\label{Bandit2ASM}
\th_{n+1}=\th_n+\g_{n+1}\left(\mathds{1}_{A_{n+1}}-\mathds{1}_{B_{n+1}}\right)h(\th_n)+\g_{n+1}\Delta M_{n+1}, \quad 
\th_0=\th\in(0,1)
\end{equation}
where $h(\th)=\th(1-\th)$, $M_n:=\sum_{k=1}^nm_k$, $M_0:=0$, with
$$m_k:=\mathds{1}_{A_{k}}(1-\th_{k-1})\left(\mathds{1}_{\{U_{n+1}\leq\th_n\}}-\th_{k-1}\right)+\mathds{1}_{B_{k}}\th_{k-1}\left((1-\th_{k-1})-\mathds{1}_{\{U_{n+1}>\th_n\}}\right).$$
We make the assumption that $A$ outperforms $B$ in average $i.e.$ that $\nu(A)>\nu(B)$ where 
\[
\frac 1n \sum_{k=1}^n \mbox{\bf 1}_{A_n}\underset{n\rightarrow\infty}{\longrightarrow} \nu(A)\quad\mbox{ and }\quad \frac 1n \sum_{k=1}^n \mbox{\bf 1}_{B_n}\underset{n\rightarrow\infty}{\longrightarrow} \nu(B)
\]
and that these convergences hold at rate $\e_n$  satisfying~(\ref{HepsASM}). Then applying Theorem \ref{Thm2ASM} with  $Y_k:=\mathds{1}_{A_{k+1}}-\mathds{1}_{B_{k+1}}$, $k\geq0$, and $\chi(y)=\frac{y}{\nu(A)-\nu(B)}$, we get a first convergence result:  as soon as $(\gamma_n)_{n\ge 1}$ is admissible in the sense of~(\ref{gammaASM})  for the sequence $(\e_n)_{n\ge1}$,
$$
\th_n\overset{a.s.}{\underset{n\rightarrow\infty}{\longrightarrow}}\th^*\in\left\{0,1\right\}.
$$
where $1$ is the target and $0$ is a trap. Further investigations on $\th^*$ are carried in \cite{TarVan} in this ergodic framework to analyze the fallibility of the algorithm which extend former results established  in~\cite{LamPagTar, LamPag2} in the purely i.i.d. setting.

\subsection{Optimal split of orders across liquidity pools}\label{DkPlASM}

This is an example of application in Finance to be implemented exclusively on  real data. It is an optimal allocation problem which solved by a stochastic   Lagrangian approach originally developed in \cite{LarLehPag}. Here, only numerical results with  real market data are presented.

\subsubsection{Model description}
The principle of a \textit{Dark pool} is to propose a price with no guarantee of executed quantity at the occasion of an OTC transaction. Usually this price is lower than the one offered on the regular market. So one can model the impact of the existence of $N$ dark pools ($N\geq 2$) on a given transaction  as follows: let $V>0$ be  the random volume to be executed, let $\th_i\in (0,1)$ be the \textit{discount factor} proposed by the dark pool $i$. Let $r_i$ denote the percentage of $V$ sent to the dark pool $i$ for execution. Let $D_i\geq 0$ be the quantity of securities that can be delivered (or made available) by the dark pool $i$ at price $\th_i S$. 

The remainder of the order is to be executed on the regular market, at price $S$. Then the  cost $C$ of the whole executed order is given by
$$C=S\sum_{i=1}^N\th_i\min\left(r_iV,D_i\right)+S\Big(V-\sum_{i=1}^N\min\left(r_iV,D_i\right)\Big) =S\Big(V-\sum_{i=1}^N\rho_i\min\big(r_iV,D_i\big)\Big)$$
where $\rho_i=1-\th_i\in(0,1), i=1,\ldots,N$. Minimizing the mean execution cost, \textit{given the price $S$}, amounts to solving the following maximization problem
\begin{equation}\label{MaxASM}
\max\Big\{\sum_{i=1}^N\rho_i\E\Big(S\min\big(r_iV,D_i\big)\Big), r\in\mathcal{P}_N\Big\}
\end{equation}
where $\mathcal{P}_N:=\Big\{r=(r_i)_{1\leq i \leq N}\in{\R}_+^N\,|\,\sum_{i=1}^Nr_i=1\Big\}$. It is then convenient to \textit{include the price $S$ into both random variables $V$ and $D_i$}  by considering $\widetilde V:=V\,S$ and $\widetilde D_i:=D_iS$ instead of $V$ and $D_i$. 

Let ${\cal I}_N=\left\{1,\ldots,N\right\}$. We set for all $r=(r_1,\ldots,r_N)\in\mathcal{P}_N$, $\Phi(r_1,\ldots,r_N):=\sum_{i=1}^N\varphi_i(r_i)$,
where
$$\forall i\in {\cal I}_N, \quad\varphi_i(u):=\rho_i\E\left(\min\left(uV,D_i\right)\right), \quad u\in\left[0,1\right].$$ 
We assume that for all $i\in{\cal I}_N$,
\begin{equation}
	V>0 \ \P\mbox{-}a.s.,\quad \P(D_i>0)>0 \, \mbox{and the distribution function of $\frac{D_i}{V}$ is continuous on ${\R}_+$,}
\label{hypVDASM}	
\end{equation}
then $\varphi_i$, $i\in {\cal I}_N$, are everywhere differentiable on the unit interval $\left[0,1\right]$ with
\begin{equation}\label{derivphiASM}
\varphi_i'(u)=\rho_i\,\E\left(\mathds{1}_{\left\{uV\leq D_i\right\}}V\right), \quad u\in\left(0,1\right],
\end{equation}
and one extends $\varphi_i$, $i\in {\cal I}_N$, on the whole real line into a concave nondecreasing function with $\lim_{\pm\infty}\varphi_i=\pm\infty$. So we can formally extend $\Phi$ on the whole affine hyperplane spanned by $\mathcal{P}_N$ $i.e.$ $\mathcal{H}_N:=\Big\{r=\left(r_1,\ldots,r_N\right)\in{\R}^N\left|\right.\sum_{i=1}^N r_i=1\Big\}$.

\subsubsection{Design of the recursive procedure}
We aim at solving the following maximization problem $\max_{r\in\mathcal{P}_N}\Phi(r)$. The Lagrangian associated to the sole affine constraint suggests that any $r^*\in\arg\max_{\mathcal{P}_N}\Phi$ iff $\varphi'_i(r^*_i)$ is constant when $i$ runs over ${\cal I}_N$ or equivalently if
$\varphi'_i(r^*_i)=\frac{1}{N}\sum_{j=1}^N\varphi'_j(r^*_j)$, $i\in {\cal I}_N$. 

We set $Y^n:=\left(V^n,D_1^n,\ldots,D_N^n\right)_{n\geq1}$. Then using the representation of the derivatives $\varphi'_i$ yields 
\begin{eqnarray*}
r^*\in\arg\max_{\mathcal{P}_N}\Phi\Longleftrightarrow \forall i\in\left\{1,\ldots,N\right\}, \ \E\Big(V\Big(\rho_i\mathds{1}_{\left\{r_i^*V<D_i\right\}}-\frac{1}{N}\sum_{j=1}^N\rho_j\mathds{1}_{\left\{r_j^*V<D_j\right\}}\Big)\Big)=0.
\end{eqnarray*}
Consequently, this leads to the following recursive zero search procedure
\begin{equation}\label{algostooptASM}
r^{n+1}_i=r^n_i+\g_{n+1}H_i(r^n,Y^{n+1}), \quad r^0\in\mathcal{P}_N, \quad n\geq0, \quad i\in {\cal I}_N,
\end{equation}
where for every $i\in {\cal I}_N$, $r\in\mathcal{P}_N$, every $V>0$ and every $D_1,\ldots,D_N\geq0$,
\begin{equation*}
H_i(r,Y)=V\Big(\rho_i\mathds{1}_{\{r_iV<D_i\}}-\frac{1}{N}\sum_{j=1}^N\rho_j\mathds{1}_{\{r_jV<D_j\}}\Big)
\end{equation*}
where $(Y^n)_{n\geq1}$ is a sequence of random vectors with nonnegative components such that, for every $n\geq1$, $(V^n,D_i^n, i=1,\ldots,N)\stackrel{d}{=}(V,D_i,i=1,\ldots,N)$.

The underlying idea of the algorithm is to reward the dark pools which outperform the mean of the $N$ dark pools by increasing the allocated volume sent at the next step (and conversely). For sake of simplicity that ${\rm argmax}_{{\cal P}_{_N}}\, \Phi= \{r^*\}\subset {\rm  int}({\cal P}_{_N})$. Our ``light'' $\nu$-averaging assumption is to assume that there exists an exponent $\eta\in (0,1]$ such that for every $u\in {\R}_+$ and every $i\!\in {\cal I}_N$
\begin{equation}\label{HypoErgoASM}
\frac{1}{n} \sum_{k=1}^n V^k  \mbox{\bf 1}_{\{u < \frac{D_i^k}{V^k}\}}-\E(V \mbox{\bf 1}_{\{u < \frac{D_i}{V}\}})= O(n^{-\eta})\quad a.s. \ \mbox{and in} \ L^2(\P)
\end{equation}
(which hold under geometric $\alpha$-mixing assumptions on $(D^n,V^n)_{n\ge1}$). Under additional technical assumptions on the support of ${\cal L}(Y^n)$ (see \cite{LarLehPag}), we can apply Theorem \ref{Thm1ASM}: if the sequence $(\g_n)_{n\geq1}$ satisfies (\ref{gammaASM}), we get that the algorithm defined by~(\ref{algostooptASM}) $a.s.$ converges towards $r^*={\rm argmax}_{{\cal P}_{_N}}\,\Phi$.

\subsubsection{Numerical Tests}
We consider the shortage setting, $i.e.$ $\E V> \sum_{i=1}^N\E D_i$ because it is the most interesting case and the most common in the market. Now, we introduce an index to measure the performances of our recursive allocation procedure.  
	
\smallskip
\noindent $\rhd$ {\bf Relative cost reduction (w.r.t. the regular market):} it is defined as the  ratios between the cost reduction of the execution using dark pools and the cost resulting from an execution on the regular market, $i.e.$, for every $n\ge 1$, 
$$\displaystyle\frac{CR^{algo}}{V^n}= \frac{\sum_{i=1}^N\rho_i\min\left(r_i^nV^n,D_i^n\right)}{V^n}.$$
We have considered for $V$ the traded volumes of a very liquid security --~namely the asset BNP~--  during an $11$ day period. Then we selected the $N$ most correlated assets (in terms of traded volumes) with the original asset. These assets are denoted $S_i$, $i=1,\ldots,N$ and we considered their traded volumes during the same 11 day period.
Finally,  the available volumes of each dark pool $i$ have been modeled  as follows using   the mixing function
$$
\forall \, 1\leq i\leq N, \quad D_i:=\beta_i\Big((1-\alpha_i)V+\alpha_i S_i \frac{\E V}{\E S_i}\Big)
$$
where $\alpha_i, \  i=1,\ldots,N$ are the recombining coefficients, $\beta_i, \  i=1,\ldots,N$ some scaling factors and $\E \,V$ and $\E\, S_i$ stand for the empirical mean of the data sets of $V$ and $S_i$. The simulations presented here have been made with four dark pools ($N=4$). Since the data used here cover $11$ days,  it is clear that, unlike the simulated data, these pseudo-real data are not stationary: in particular they are subject to daily changes of trend and volatility (at least). To highlight the resulting changes in the response of the algorithms, we have specified the days by drawing vertical doted lines. The dark pool pseudo-data parameters are set to $\beta=(0.1,0.2,0.3,0.2)^t$, $\alpha=(0.4,0.6,0.8,0.2)^t$ and the dark pool trading (rebate) parameters are set to $\rho=(0.0 , 0.02, 0.04, 0.06)^t$.

\begin{figure}[!ht]
\centering
\includegraphics[width=10cm]{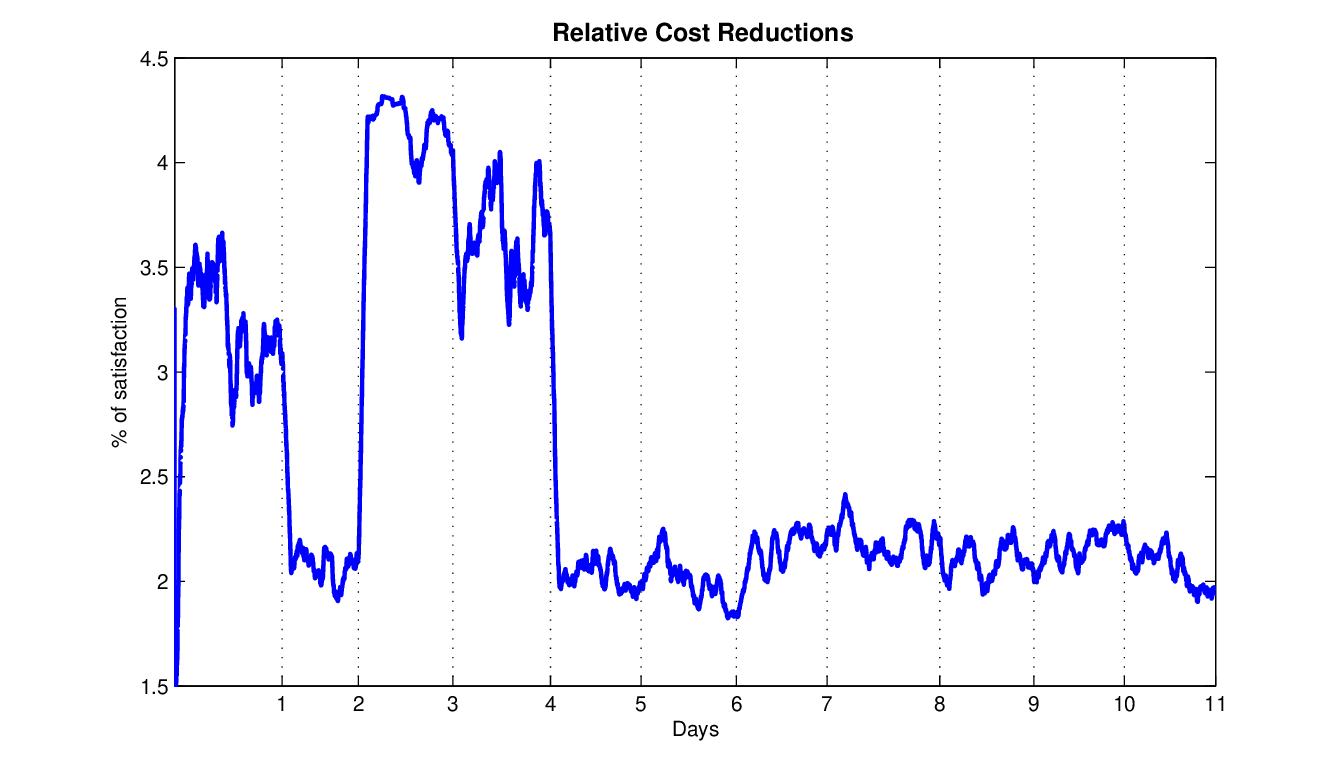}
\vspace{-0.5cm}
\caption{Case $N=4$, $\sum_{i=1}^N\beta_i<1$, $0.2<\alpha_i\leq0.8$ and $r^0_i=1/N$, $1\leq i\leq N$.}
\label{fig:3}
\end{figure}

We benchmarked  --~see Figure~\ref{fig:3}~-- the algorithm on the whole data set ($11$ days) as though it were stationary. In particular, the running means of the performances are computed from the very beginning for the first 1500 data, and then by a moving average computed on a window of 1500 data. 

\vspace{-0.25cm}

\vspace{-0.5cm}

\small
\bibliographystyle{plain}
\bibliography{biblio}

\begin{thebibliography}{10}

\bibitem{Ben}
M.~Ben Alaya.
\newblock Sur la m\'ethode du shift en simulation.
\newblock In {\em Probabilit\'es num\'eriques}, volume~10 of {\em Collect.
  Didact.}, pages 61--66. INRIA, Rocquencourt, 1992.

\bibitem{Ben1}
M.~Ben Alaya.
\newblock On the simulation of expectations of random variables depending on a
  stopping time.
\newblock {\em Stochastic Anal. Appl.}, 11(2):133--153, 1993.

\bibitem{BenPag}
M.~Ben Alaya and G.~Pag{\`e}s.
\newblock Rate of convergence for computing expectations of stopping
  functionals of an {$\alpha$}-mixing process.
\newblock {\em Adv. in Appl. Probab.}, 30(2):425--448, 1998.

\bibitem{AndMouPri}
C.~Andrieu, \'E. Moulines, and P.~Priouret.
\newblock Stability of stochastic approximation under verifiable conditions.
\newblock {\em SIAM J. Control Optim.}, 44(1):283--312, 2005.

\bibitem{BarFriPag}
O.~Bardou, N.~Frikha, and G.~Pag{\`e}s.
\newblock Computing {V}a{R} and {CV}a{R} using stochastic approximation and
  adaptive unconstrained importance sampling.
\newblock {\em Monte Carlo Methods Appl.}, 15(3):173--210, 2009.

\bibitem{BMP}
A.~Benveniste, M.~M{\'e}tivier, and P.~Priouret.
\newblock {\em Adaptive algorithms and stochastic approximations}, volume~22 of
  {\em Applications of Mathematics (New York)}.
\newblock Springer-Verlag, Berlin, 1990.
\newblock Translated from the French by Stephen S. Wilson.

\bibitem{BoLe}
N.~Bouleau and D.~L\'epingle.
\newblock {\em Numerical methods for stochastic processes}.
\newblock Wiley Series in Probability and Mathematical Statistics. John Wiley
  \& Sons Inc., New York, first edition, 1994.
\newblock A Wiley-Interscience Publication.

\bibitem{Bru}
M.-F. Bru.
\newblock Wishart processes.
\newblock {\em J. Theoret. Probab.}, 4(4):725--751, 1991.

\bibitem{DedDou}
J.~Dedecker, P.~Doukhan, G.~Lang, J.~R. Le{\'o}n~R., S.~Louhichi, and
  C.~Prieur.
\newblock {\em Weak dependence: with examples and applications}, volume 190 of
  {\em Lecture Notes in Statistics}.
\newblock Springer, New York, 2007.

\bibitem{DedMerVol}
J.~Dedecker, F.~Merlev{\`e}de, and D.~Voln{\'y}.
\newblock On the weak invariance principle for non-adapted sequences under
  projective criteria.
\newblock {\em J. Theoret. Probab.}, 20(4):971--1004, 2007.

\bibitem{DFMP}
R.~Douc, G.~Fort, \'E. Moulines, and P.~Priouret.
\newblock Forgetting the initial distribution for hidden {M}arkov models.
\newblock {\em Stochastic Process. Appl.}, 119(4):1235--1256, 2009.

\bibitem{Dou}
P.~Doukhan.
\newblock {\em Mixing}, volume~85 of {\em Lecture Notes in Statistics}.
\newblock Springer-Verlag, New York, 1994.
\newblock Properties and examples.

\bibitem{Duf}
M.~Duflo.
\newblock {\em Algorithmes stochastiques}, volume~23 of {\em Math\'ematiques \&
  Applications (Berlin) [Mathematics \& Applications]}.
\newblock Springer-Verlag, Berlin, 1996.

\bibitem{FolSch}
H.~F{\"o}llmer and A.~Schied.
\newblock {\em Stochastic finance}, volume~27 of {\em de Gruyter Studies in
  Mathematics}.
\newblock Walter de Gruyter \& Co., Berlin, 2002.
\newblock An introduction in discrete time.

\bibitem{Fri}
N.~Frikha.
\newblock {\em Contribution \`a la mod\'elisation et \`a la gestion dynamique
  du risque des march\'es de l'\'energie}.
\newblock PhD thesis, UPMC, 2010.

\bibitem{GalKok}
I.~S. G{\'a}l and J.~F. Koksma.
\newblock Sur l'ordre de grandeur des fonctions sommables.
\newblock {\em C. R. Acad. Sci. Paris}, 227:1321--1323, 1948.

\bibitem{GouJasSuf}
C.~Gouri{\'e}roux, J.~Jasiak, and R.~Sufana.
\newblock The {W}ishart autoregressive process of multivariate stochastic
  volatility.
\newblock {\em J. Econometrics}, 150(2):167--181, 2009.

\bibitem{KusCla}
H.~J. Kushner and D.~S. Clark.
\newblock {\em Stochastic approximation methods for constrained and
  unconstrained systems}, volume~26 of {\em Applied Mathematical Sciences}.
\newblock Springer-Verlag, New York, 1978.

\bibitem{KusYin}
H.~J. Kushner and G.~G. Yin.
\newblock {\em Stochastic approximation and recursive algorithms and
  applications}, volume~35 of {\em Applications of Mathematics (New York)}.
\newblock Springer-Verlag, New York, second edition, 2003.
\newblock Stochastic Modelling and Applied Probability.

\bibitem{LamPag}
D.~Lamberton and G.~Pag{\`e}s.
\newblock Recursive computation of the invariant distribution of a diffusion.
\newblock {\em Bernoulli}, 8(3):367--405, 2002.

\bibitem{LamPag2}
D.~Lamberton and G.~Pag{\`e}s.
\newblock How fast is the bandit?
\newblock {\em Stoch. Anal. Appl.}, 26(3):603--623, 2008.

\bibitem{LamPag3}
D.~Lamberton and G.~Pag{\`e}s.
\newblock A penalized bandit algorithm.
\newblock {\em Electron. J. Probab.}, 13:no. 13, 341--373, 2008.

\bibitem{LamPagTar}
D.~Lamberton, G.~Pag{\`e}s, and P.~Tarr{\`e}s.
\newblock When can the two-armed bandit algorithm be trusted?
\newblock {\em Ann. Appl. Probab.}, 14(3):1424--1454, 2004.

\bibitem{LapPagSab}
B.~Lapeyre, G.~Pag{\`e}s, and K.~Sab.
\newblock Sequences with low discrepancy---generalisation and application to
  {R}obbins-{M}onro.
\newblock {\em Statistics}, 21(2):251--272, 1990.

\bibitem{Lar}
S.~Laruelle.
\newblock {\em Analyse d'algorithmes stochastiques appliqu\'es \`a la Finance}.
\newblock PhD thesis, UPMC, 2011.

\bibitem{LarLehPag}
S.~Laruelle, C.-A. Lehalle, and G.~Pag{\`e}s.
\newblock Optimal split of orders across liquidity pools: a stochastic
  algorithm approach.
\newblock Pre-pub LPMA 1316, to appear in {\em SIAM J. of Financial
  Mathematics}, 2009.

\bibitem{Lem}
V.~Lemaire.
\newblock {\em Estimation r\'ecursive de la mesure invariante d'un processus de
  diffusion}.
\newblock PhD thesis, Universit\'e de Marne-La-Vall\'ee, 2005.

\bibitem{Lem2}
V.~Lemaire.
\newblock An adaptive scheme for the approximation of dissipative systems.
\newblock {\em Stochastic Processes and their Applications},
  117(10):1491--1518, 2007.

\bibitem{LemPag}
V.~Lemaire and G.~Pag{\`e}s.
\newblock Unconstrained recursive importance sampling.
\newblock {\em Ann. Appl. Probab.}, 20(3):1029--1067, 2010.

\bibitem{LokZer}
A.~L{\o}kka and M.~Zervos.
\newblock A model for the long-term optimal capacity level of an investment
  project.
\newblock {\em Int. J. Theor. Appl. Finance}, 14(2):187--196, 2011.

\bibitem{MehMey}
P.~G. Mehta and S.~P. Meyn.
\newblock Q-learning and {Pontryagin's} minimum principle.
\newblock In {\em Proceedings of the 48th IEEE Conference on Decision and
  Control. Held jointly with the 2009 28th Chinese Control Conference. CDC/CCC
  2009.}, pages 3598--3605, Dec. 2009.

\bibitem{NarSha}
K.~S. Narendra and I.~J. Shapiro.
\newblock Use of stochastic automata for parameter self-optimization with
  multimodal performance criteria.
\newblock {\em IEEE Trans. Syst. Sci. and Cybern.}, 5(4):352 -- 360, 1969.

\bibitem{Nie}
H.~Niederreiter.
\newblock {\em Random number generation and quasi-{M}onte {C}arlo methods},
  volume~63 of {\em CBMS-NSF Regional Conference Series in Applied
  Mathematics}.
\newblock Society for Industrial and Applied Mathematics (SIAM), Philadelphia,
  PA, 1992.

\bibitem{Nor}
M.~F. Norman.
\newblock On the linear model with two absorbing barriers.
\newblock {\em J. Mathematical Psychology}, 5:225--241, 1968.

\bibitem{PelUte}
M.~Peligrad and S.~Utev.
\newblock Central limit theorem for stationary linear processes.
\newblock {\em Ann. Probab.}, 34(4):1608--1622, 2006.

\bibitem{Pro}
P.~D. Pro{\u\i}nov.
\newblock Discrepancy and integration of continuous functions.
\newblock {\em J. Approx. Theory}, 52(2):121--131, 1988.

\bibitem{RobMon}
H.~Robbins and S.~Monro.
\newblock A stochastic approximation method.
\newblock {\em Ann. Math. Statistics}, 22:400--407, 1951.

\bibitem{RocUry}
R.~T. Rockafellar and S.~Uryasev.
\newblock Conditional value-at-risk: optimization approach.
\newblock In {\em Stochastic optimization: algorithms and applications
  ({G}ainesville, {FL}, 2000)}, volume~54 of {\em Appl. Optim.}, pages
  411--435. Kluwer Acad. Publ., Dordrecht, 2001.

\bibitem{ShiMey}
D.~Shirodkar and S.~Meyn.
\newblock Quasi stochastic approximation.
\newblock In {\em American Control Conference, 2011. ACC '11.}, June 2011.

\bibitem{TarVan}
P.~Tarr{\`e}s and P.~Vandekerkhove.
\newblock On ergodic two-armed bandits.
\newblock To appear in {\em Annals of Applied Probability}, 2009.

\end{thebibliography}

\end{document}